\pdfoutput=1
\RequirePackage{ifpdf}
\ifpdf 
\documentclass[pdftex]{sigma}
\else
\documentclass{sigma}
\fi

\numberwithin{equation}{section}
\usepackage{mathtools}
\usepackage{tikz-cd}

\newtheorem{Theorem}{Theorem}[section]
\newtheorem*{Theorem*}{Theorem}
\newtheorem{Corollary}[Theorem]{Corollary}
\newtheorem{Lemma}[Theorem]{Lemma}
\newtheorem{Proposition}[Theorem]{Proposition}

\theoremstyle{definition}

\begin{document}

\allowdisplaybreaks

\newcommand{\arXivNumber}{2505.02954}

\renewcommand{\thefootnote}{}

\renewcommand{\PaperNumber}{048}

\FirstPageHeading

\ShortArticleName{1-Point Functions for $\mathbb{Z}_2$-Orbifolds of Lattice VOAs}

\ArticleName{1-Point Functions for $\boldsymbol{\mathbb{Z}_2}$-Orbifolds of Lattice VOAs\footnote{This paper is a~contribution to the Special Issue on Recent Advances in Vertex Operator Algebras in honor of James Lepowsky. The~full collection is available at \href{https://sigma-journal.com/Lepowsky.html}{https://sigma-journal.com/Lepowsky.html}}}

\Author{Maneesha AMPAGOUNI}

\AuthorNameForHeading{M.~Ampagouni}

\Address{Department of Mathematics, University of California, Santa Cruz, \\
1156 High St, Santa Cruz, CA 95064, USA}
\Email{\mail{dampagou@ucsc.edu}}
\URLaddress{\url{https://sites.google.com/view/maneeshaampagouni}}

\ArticleDates{Received May 07, 2025, in final form May 03, 2026; Published online May 13, 2026}

\Abstract{In this paper, we compute the 1-point correlation functions of all states for the $\mathbb{Z}_2$-orbifolds of lattice vertex operator algebras.}


\Classification{17B69; 11F11; 11F27}

\renewcommand{\thefootnote}{\arabic{footnote}}
\setcounter{footnote}{0}

\section{Introduction}

The first example of an orbifold in the theory of vertex operator algebras is the moonshine module $V^{\natural}$. It was constructed by Igor Frenkel, James Lepowsky and Arne Meurman \cite{Frenkel1988VertexMonster} as a~$\mathbb{Z}_2$-orbifold of the vertex operator algebra associated with the Leech lattice. The construction of the moonshine module can be generalized by replacing the Leech lattice with any even, positive-definite lattice $L$ of rank~$k=8l$, where $l\in \mathbb{Z}^+$ (the set of positive integers) such that $\sqrt{2}L^*$ is also even, where $L^*$ denotes the dual lattice of $L$, to obtain the $\mathbb{Z}_2$-orbifold of the lattice vertex operator algebra $V_L$ with respect to the involutive automorphism~$\theta$, where~$\theta$ is the lift of the $(-1)$-involution of the lattice~$L$. This was first done explicitly by Dolan, Goddard and Montague~\cite{Dolan1990ConformalOperators}. This condition on $L$ is more general than unimodularity. When $L$ is unimodular, the orbifold that is obtained through this construction indeed coincides with that of the holomorphic vertex operator algebra obtained through the cyclic orbifold theory given by Ekeren, M{\"o}ller and Scheithauer in~\cite{van2020Construction}.

One of the most intriguing features of the moonshine module as it was initially constructed was its character (1-point function corresponding to the vacuum state) being the modular function $j(\tau)-744$. The modular invariance of the character function was later explained by Zhu~\cite{Zhu1996ModularAlgebras} as a consequence of the axioms of the vertex operator algebras together with some finiteness conditions. Thus modular invariance of characters of vertex operator algebras (VOAs) became a subject of interest further. In \cite{dong2001Quasi}, Dong, Mason and Nagatomo investigated the modular properties of 1-point functions corresponding to free bosonic VOAs and lattice VOAs and found that the trace functions in these two theories have the shape $f(q)/\eta(q)^d$, where $f(q)$ is quasi-modular in the case of $d$ free bosons and modular in the latter case. In analogy with these works, one may study the modular properties of 1-point functions associated with the $\mathbb{Z}_2$-orbifolds of lattice VOAs. This requires, as a first step, explicit computations of such 1-point functions, which is the focus of the present paper.

In this paper, we focus on the computation of 1-point correlation functions for the \( \mathbb{Z}_2 \)-orbifolds of vertex operator algebras associated with unimodular, even, positive definite lattices $L$ of rank~${k=8l}$ where $l\in \mathbb{Z}^+$. The approach in this paper is inspired by the techniques developed by Geoffrey Mason and Michael Mertens in their paper on 1-point functions for symmetrized Heisenberg and lattice vertex operator algebras \cite{Mason20231-PointAlgebras}. We adopt and extend those techniques to the \( \mathbb{Z}_2 \)-orbifold setting to perform computations on the twisted sector and obtain results for the 1-point functions for the $\mathbb{Z}_2$-orbifold. While the work by Mason and Mertens sheds light on how the 1-point functions traced over $V_L^+$ exhibit modular invariance with respect to a congruence subgroup, our work finds the trace contributed by \smash{$\bigl(V_L^T\bigr)^+$} (which is the $+1$-eigenspace of the action of the $(-1)$-involution ($\theta$) on $V_L^T$ defined in Section~\ref{sec2.3}). Combining the two traces, we obtain the 1-point functions for the $\mathbb{Z}_2$-orbifold~$V$~\cite{van2020Construction}, which are expected to exhibit modular invariance under the full modular group up to a character. This modular invariance is a property that can be realized as a special case of \cite[Theorem 5.1.1]{Zhu1996ModularAlgebras} when~${n=1}$ and it is verified in Proposition~\ref{proposition1}. Further, we also observe that due to the structure of the twisted module, the 1-point functions corresponding to the lattice states of the form $h_{i_1}[-n_1]\cdots h_{i_p}[-n_p]e_{\alpha}$ when~${\alpha\in L\backslash 2L}$ vanish.

The main contributions of this paper are as follows:

1. We derive explicit formulas for the 1-point correlation functions of the \( \mathbb{Z}_2 \)-orbifold of lattice vertex operator algebras, as in the following theorems.
 \begin{Theorem}
Let $L$ be a positive-definite even unimodular lattice of rank $k=8l$ {\rm(}where $l\in \mathbb{Z}^+${\rm)}. Let $V$ be the $\mathbb{Z}_2$-orbifold of the VOA $V_L$ associated with the lattice $L$ formed by the construction of Dolan, Goddard and Montague {\rm\cite{Dolan1990ConformalOperators} (}also by {\rm\cite{van2020Construction})}.
For positive integers $n_i\geq 1$, corresponding to the Heisenberg state
$u\coloneqq h_{i_1}[-n_1]\cdots h_{i_p}[-n_p]\mathbf{1}$,
we have the $1$-point function given by
\begin{align*}
Z_V(u,\tau)={}&\frac{1}{2}\sum_{\Delta\subseteq \Lambda}\frac{\theta_L(\tau,P_{\Delta})}{\eta(\tau)^k}\bigg(\sum_{\sigma\in \operatorname{Inv}_0(\underline{p}\backslash \Delta)}\prod_{(rs)}\delta_{i_r,i_s}\hat{E}_{n_r+n_s}(\tau)\bigg)\\
&+\frac{1}{2}\eta(\tau)^{k/2}\left(\frac{\Theta_1(\tau)}{2}\right)^{-k/2}\bigg(\sum_{\sigma\in \operatorname{Inv}_0(\underline{p})} \prod_{(rs)}\delta_{i_r,i_s}\hat{F}_{n_r+n_s}(\tau) \bigg)\\
&+\frac{1}{2}\eta(\tau)^{k/2}\left(\frac{\Theta_2(\tau)}{2}\right)^{-k/2}\bigg(\sum_{\sigma\in \operatorname{Inv}_0(\underline{p})}\prod_{(rs)}\delta_{i_r,i_s}\overline{E}_{n_r+n_s}(\tau)\bigg)\\
&+\frac{(-1)^l}{2}\eta(\tau)^{k/2}\left(\frac{\Theta_3(\tau)}{2}\right)^{-k/2}\bigg(\sum_{\sigma\in \operatorname{Inv}_0(\underline{p})}\prod_{(rs)}\delta_{i_r,i_s}\overline{F}_{n_r+n_s}(\tau)\bigg),
\end{align*}
where $\Lambda=\{j\in \underline{p}\mid n_j=1\}$, $\Delta$ is a subset of $\Lambda$ of even cardinality, $\underline{p}\coloneqq \{1,2,\dots, p\}$, for a set~$S$, $\operatorname{Inv}_0(S)$ is the set of all fixed-point-free involutions of the set $S$ and
$P_{\Delta}(\alpha)=\prod_{j\in\Delta}\langle h_{i_j},\alpha\rangle$.
More details about the notation can be found in Section {\rm\ref{section2}}.
 \end{Theorem}

 \begin{Theorem}
In the same setting as the above theorem, for positive integers $n_i\geq 1$, $\alpha\in 2L$, corresponding to the state
$u\coloneqq h_{i_1}[-n_1]\cdots h_{i_p}[-n_p]e_{\alpha}$,
where $e_{\alpha}$ is $f_{\alpha}={\rm e}^{\alpha}+{\rm e}^{-\alpha}$ or $g_{\alpha}={\rm e}^{\alpha}-{\rm e}^{-\alpha}$ accordingly as $p$ is even or odd, respectively, the $1$-point function is given by
\begin{align*}
Z_V(u,\tau)={}&\eta(\tau)^{k/2}\left(\frac{\Theta_1(\tau)}{2}\right)^{\langle\alpha,\alpha\rangle-k/2}\bigg(\sum_{\sigma} \prod_{(rs)(t)}\delta_{i_r,i_s}\Tilde{F}_{n_t}(\tau)\hat{F}_{n_r+n_s}(\tau)\bigg)\\
&+\eta(\tau)^{k/2}\left(\frac{\Theta_2(\tau)}{2}\right)^{\langle\alpha,\alpha\rangle-k/2}
\bigg(\sum_{\sigma}\prod_{(rs)(t)}\delta_{i_r,i_s}\Tilde{\overline{E}}_{n_t}\overline{E}_{n_r+n_s}(\tau)\bigg)\\
&+(-1)^l\eta(\tau)^{k/2}\left(\frac{\Theta_3(\tau)}{2}\right)^{\langle\alpha,\alpha\rangle-k/2}
\bigg(\sum_{\sigma}\prod_{(rs)(t)}\delta_{i_r,i_s}\Tilde{\overline{F}}_{n_t}\overline{F}_{n_r+n_s}(\tau)\bigg),
\end{align*}
where $\sigma$ ranges over all involutions of the index set $S_p=\{1,2,\dots, p\}$ and $(rs)$, $(t)$ range over the $2$-cycles and $1$-cycles respectively in the decomposition of $\sigma$ in $S_p$ as a product of disjoint $2$-cycles and $1$-cycles. For $\alpha\in L\backslash 2L$, we have
$Z_V(u,\tau)=0$.
\end{Theorem}

2. We analyse the modular properties of these 1-point functions, confirming that they indeed exhibit modular invariance under the full modular group (up to a character).

This paper is organized as follows: In Section~\ref{section2}, we provide a brief overview of the necessary background and notation on modular forms, elliptic-type functions, theta functions, lattice vertex operator algebras, the corresponding $\mathbb{Z}_2$-orbifold, and the square bracket formalism. Section~\ref{sec3} outlines the methods used by Mason and Mertens, which we adapt for our purposes to develop $\mathbb{Z}_2$-twisted Zhu theory. Here, we produce a reduction theorem for 2-point functions traced over the twisted module $V_L^T$ and similar formulas were first produced in \cite[Theorem~8.4]{Dong2000Modular}. In Section~\ref{sec4}, we present some computations in the twisted space. In Section~\ref{sec5}, we present our main results, including the explicit computation of 1-point functions for the \( \mathbb{Z}_2 \)-orbifolds of lattice vertex operator algebras. In Section~\ref{sec6}, we verify that these 1-point functions exhibit modular invariance under the full modular group up to a character.

\section{Background and notation}\label{section2}

\subsection{Modular forms and elliptic-type functions}

Let $\mathbb{H}$ denote the complex upper-half plane, $\tau$ a typical element in $\mathbb{H}$ and $q={\rm e}^{2\pi {\rm i} \tau}$. The Dedekind eta function is defined as
\[\eta(\tau)=q^{1/24}\prod_{n\geq 1}(1-q^n).\]
The Bernoulli numbers $B_k$ are defined by
\[\sum_{k\geq 0}\frac{B_k}{k!}z^k\coloneqq \frac{z}{({\rm e}^z-1)}.\]
The Eisenstein series for even positive integers $k$ are defined by
\[E_k(\tau)\coloneqq -\frac{B_k}{k!}+\frac{2}{(k-1)!}\sum_{n\geq 1}\sigma_{k-1}(n)q^n,\]
where $\sigma_{k-1}(m)\coloneqq \sum_{d\mid m}d^{k-1}$ denotes the usual divisor sum function. We denote the level 2 Eisenstein series defined for positive $k$, by $F_k$, where
$F_k(\tau)\coloneqq 2E_k(2\tau)-E_k(\tau)$.
Note that~$E_k(\tau)$ is a modular form of weight $k$ for the full modular group ${\rm SL}_2(\mathbb{Z})$ as long as~${k\neq 2}$ and the functions $F_k$ are modular forms of weight $k$ for
\[\Gamma_0(2) \coloneqq \left\{\begin{pmatrix}
 a & b \\
 c & d
\end{pmatrix}\mid c\equiv 0 \ (\text{mod }2)\right\}.\]
We use the same notation for renormalised Eisenstein series as in \cite{Mason20231-PointAlgebras}
\begin{gather*}
\hat{E}_{m+n}(\tau)\coloneqq (-1)^{n+1}n \binom{m+n-1}{n} E_{m+n}(\tau),\\
\hat{F}_{m+n}(\tau)\coloneqq (-1)^{n+1}n \binom{m+n-1}{n} F_{m+n}(\tau).
\end{gather*}
Additionally, we introduce a new notation for the following level 2 and level 4 modular forms
\begin{gather*}
\overline{E}_{m+n}(\tau)\coloneqq 2^{-(m+n-1)}\hat{E}_{m+n}(\tau/2)-\hat{E}_{m+n}(\tau),\\
\overline{F}_{m+n}(\tau)\coloneqq 2^{-(m+n-1)}\hat{F}_{m+n}(\tau/2)-\hat{F}_{m+n}(\tau).
\end{gather*}
We also use the same notation for the following elliptic-type functions used in \cite{Mason20231-PointAlgebras} as
\begin{equation}\label{eq 7 from MM}
 P_1(z,\tau)\coloneqq \frac{1}{2}+\sum_{0\neq n\in \mathbb{Z}}\frac{{\rm e}^{nz}}{(1-q^n)}=-\frac{1}{z}+\sum_{k\geq 1}E_{2k}(\tau)z^{2k-1}=E_{2}(\tau)z-\zeta(z,\tau),
\end{equation}
where $\zeta(z,\tau)$ is the Weirstrass zeta function for the lattice $\Lambda=2\pi{\rm i}(\mathbb{Z}\oplus\mathbb{Z}\tau)$, i.e.,
\[\zeta(z,\tau)=\frac{1}{z}+\sum_{\omega\in \Lambda}\left(\frac{1}{z-\omega}+\frac{1}{\omega}+\frac{z}{\omega^2} \right),\qquad z\in \mathbb{C}\backslash \Lambda.\]
The $z$-derivatives of the above Lambert series are thus
\[P_1^{(m)}(z,\tau)=m!\left( \frac{(-1)^{m+1}}{z^{m+1}}+\sum_{k\geq 1}\binom{2k-1}{m}E_{2k}(\tau)z^{2k-m-1} \right),\qquad m\geq 0.\]
We also use the function
\[Q_1(z,\tau)\coloneqq \sum_{n\in \mathbb{Z}}\frac{{\rm e}^{nz}}{1+q^n},\]
whose higher $z$-derivatives have series expansions as follows:
\begin{gather}
 Q_1^{(m)}(z,\tau)\coloneqq \sum_{n\in \mathbb{Z}}\frac{n^m{\rm e}^{nz}}{1+q^n}=m!\left( \frac{(-1)^{m+1}}{z^{m+1}}+\sum_{k\geq 1}\binom{2k-1}{m}F_{2k}(\tau)z^{2k-m-1} \right)\label{eq 10 from MM}
\end{gather}
 for $m\geq 0$.

\subsection{Theta functions}

The theta function of an even lattice $L$ of rank $k$ with a positive definite bilinear form $\langle \cdot,\cdot\rangle$ is defined by
\[\theta_L(\tau)=\sum_{\alpha\in L}q^{\langle\alpha,\alpha\rangle/2}.\]
We use the same notation as \cite{Mason20231-PointAlgebras} for the following variant of the theta function:
\[\theta_L(\tau,v,m)=\sum_{\alpha\in L}\langle v,\alpha\rangle^mq^{\langle\alpha,\alpha\rangle/2}\]
for any vector $v\in L\otimes\mathbb{C}$ and non-negative integer $m\geq 0$. We also require the following variant of the theta function of a lattice $L$. For any function $P\colon L\to \mathbb{C}$, set
\[\theta_L(\tau,P)\coloneqq \sum_{\alpha\in L}P(\alpha)q^{\langle\alpha,\alpha\rangle/2}.\]
We follow the notation of \cite{Dong2000MonstrousWeight} to denote the Jacobi theta functions by $\Theta_1$, $\Theta_2$, $\Theta_3$, where
\[\Theta_1(\tau)=2\frac{\eta(2\tau)^2}{\eta(\tau)}, \qquad \Theta_2(\tau)=\frac{\eta(\tau/2)^2}{\eta(\tau)},\qquad \Theta_3(\tau)=\frac{\eta(\tau)^5}{\eta(\tau/2)^2\eta(2\tau)^2}.\]

\subsection[Lattice vertex operator algebra V\_L and the Z\_2-twisted module V\_L^T]{Lattice vertex operator algebra $\boldsymbol{V_L}$ and the $\boldsymbol{\mathbb{Z}_2}$-twisted module $\boldsymbol{V_L^T}$}\label{sec2.3}

Let $L$ be an even, positive definite, unimodular lattice of rank $k=8l$ (for some $l\in \mathbb{Z}^+$) with the bilinear form
$\langle\cdot,\cdot \rangle\colon L\times L\to \mathbb{Z}$. The $\mathbb{C}$-linear extension defines a bilinear form on $\mathfrak{h}\coloneqq L\otimes \mathbb{C}$. Let $\bigl(\hat{L},^-\bigr)$ be a central extension of $L$ by a cyclic group of order 2 $\langle \kappa\rangle=\big\langle \kappa \mid \kappa^2=1\big\rangle$
such that
\[\begin{tikzcd}
 1 \arrow[r] & \langle\kappa\rangle \arrow[r] & \hat{L} \arrow[r, "-"] & L \arrow[r] & 1.
\end{tikzcd}\]
Let
$e\colon L\to \hat{L}$, $ \alpha\to e_{\alpha}$
be a section of $\hat{L}$, that is, a map $e$ such that $^-\circ e=1$ such that
$e_0=1$,
and denote by $\epsilon_0\colon L\times L\to \mathbb{Z}/2\mathbb{Z}$
the corresponding $2$-cocycle, which is defined by the condition
\[e_{\alpha}e_{\beta}=\kappa^{\epsilon_0(\alpha,\beta)}e_{\alpha+\beta}\qquad \text{for} \ \alpha,\beta\in L.\]
Further, define
\[\epsilon\colon\ L\times L\to \mathbb{C}^{\times},\qquad
(\alpha,\beta)\to (-1)^{\epsilon_0(\alpha,\beta)}.\]
We can choose $\epsilon$ \cite{Mason2018N2N4Subalgebras} (and hence $\epsilon_0$) so that it satisfies
\[\epsilon(\alpha,\beta)\epsilon(\beta,\alpha)=(-1)^{\langle\alpha,\beta\rangle+\langle\alpha,\alpha\rangle\langle\beta,\beta\rangle},\qquad \epsilon(\alpha,\alpha)=(-1)^{(\langle\alpha,\alpha\rangle+\langle\alpha,\alpha\rangle^2)/2}.\]
Define a faithful character
$\chi\colon \langle\kappa\rangle\to \mathbb{C}^{\times}$
by the condition
$\chi(\kappa)=-1$.
Denote by $\mathbb{C}_{\chi}$ the one-dimensional space $\mathbb{C}$ viewed as a $\langle\kappa\rangle$-module on which $\langle\kappa\rangle$ acts according to $\chi$:
$\kappa.1=-1$.
Denote by \smash{$\mathbb{C}\{L\}={\rm Ind}_{\langle\kappa\rangle}^{\hat{L}}\mathbb{C}_{\chi} =\mathbb{C}[\hat{L}]\otimes_{\mathbb{C}[\langle\kappa\rangle]}\mathbb{C}_{\chi}
=\mathbb{C}[\hat{L}]/(\kappa+1)\mathbb{C}[\hat{L}]$}. For $a\in \hat{L}$, set
$\iota(a)=a\otimes 1\in \mathbb{C}\{L\}$.
The choice of the section allows us to identify $\mathbb{C}\{L\}$ with the group algebra $\mathbb{C}[L]$ twisted by the 2-cocycle $\epsilon$, viewed as a vector space, by the linear isomorphism \cite{Lepowsky2004Vertex}
\[\mathbb{C}[L]\to \mathbb{C}\{L\},\qquad {\rm e}^{\alpha}\to \iota(e_{\alpha})\qquad \text{for}\ \alpha\in L,\]
where ${\rm e}^{\alpha}$ is the image of $1$ through $e_{\alpha}$ (realised through the action of $\hat{L}$ on $\mathbb{C}[L]$). Recall that for $\alpha,\beta, \gamma\in L$, the action of $\hat{L}$ on $\mathbb{C}[L]$ is given by
\[e_{\alpha}.{\rm e}^{\beta}=\epsilon(\alpha,\beta){\rm e}^{\alpha+\beta},\qquad\kappa. {\rm e}^{\beta}=-{\rm e}^{\beta}\qquad \text{for} \ \alpha,\beta\in L.\]
Thus following the notation as in \cite{Mason20231-PointAlgebras}, we denote the Fock space of the corresponding Heisenberg VOA of rank $k$ with $M$ and that of the corresponding lattice theory VOA with
\[V_L= M\otimes \mathbb{C}\{L\}= M\otimes \mathbb{C}[L]=\oplus_{\alpha\in L}M\otimes {\rm e}^{\alpha}.\]
We further use the same notation as in \cite[Section 2.4]{Mason20231-PointAlgebras} for representing different Fock states of~$V_L$ after considering an orthonormal basis $\{h_1,h_2,\dots,h_k\}$ of $\mathfrak{h}$. However, we use $\theta$ to denote the involutive automorphism of $V_L$ which is a lift of the negating automorphism of $L$. To form the $\mathbb{Z}_2$-twisted module $V_L^T$, consider the twisted affine algebra
\[\hat{\mathfrak{h}}[-1]=\bigoplus_{n\in \mathbb{Z}+1/2}\mathfrak{h}\otimes t^n \oplus \mathbb{C}c.\]
Let $S\bigl(\hat{\mathfrak{h}}[-1]^-\bigr)$ denote the symmetric algebra generated by the $\mathbb{Z}^-$-graded subalgebra of the twisted affine algebra. Let the objects associated with \smash{$\hat{L}$} be defined as in \cite[equations (7.1.6)--(7.1.27)]{Frenkel1988VertexMonster} and let $T$ be any irreducible $\hat{L}$-module such that
$\kappa.v=-v$ for $v\in T$.
Observe that we set $s=2$ and thus $\omega=-1$ in the notation of \cite{Frenkel1988VertexMonster}. Since $L$ is unimodular, there exists a unique such $T$. The $\mathbb{Z}_2$-twisted module $V_L^T$ is the space
\[V_L^T=S(\mathfrak{h}[-1]^-)\otimes T.\]
Note that the elements of $T$ are all graded by the $L[0]$-weight $k/16$ \cite{Frenkel1988VertexMonster}.
Further, for $h_i\in \mathfrak{h}$, $n_i-1/2\in \mathbb{Z}^+\cup\{0\}$, $t\in T$, the action of $\theta$ on the twisted module $V_L^T$ is defined by
\[\theta(h_{i_1}[-n_1]\cdots h_{i_j}[-n_j]\otimes t)=(-1)^{j+k/8} h_{i_1}[-n_1]\cdots h_{i_j}[-n_j]\otimes t.\]

\subsection[Z\_2-orbifold of the lattice VOA V\_L]{$\boldsymbol{\mathbb{Z}_2}$-orbifold of the lattice VOA $\boldsymbol{V_L}$}

The holomorphic $\mathbb{Z}_2$-orbifold $V$ of $V_L$ is thus formed from the fixed subspaces of the involutive automorphism $\theta$, $V_L^{\theta}\subseteq V_L$ and $\bigl(V_L^T\bigr)^{\theta}\subseteq V_L^T$ as \cite{Dolan1990ConformalOperators, van2020Construction}
\smash{$
V=V_L^{\theta}\oplus \bigl(V_L^T\bigr)^{\theta}$}.
Moving forward, we shall denote \smash{$V_L^{\theta}$} with \smash{$V_L^+$} and \smash{$\bigl(V_L^T\bigr)^{\theta}$} with \smash{$\bigl(V_L^T\bigr)^+$} and the $-1$-eigenspaces of $\theta$ in $V_L$ and \smash{$V_L^T$} with $V_L^-$ and \smash{$\bigl(V_L^T\bigr)^-$}, respectively. We further use the same notation as in \cite{Mason20231-PointAlgebras} for states in \smash{$V_L^+$}, where for $\alpha\in L$
$
f_{\alpha}\coloneqq {\rm e}^{\alpha}+{\rm e}^{-\alpha}$, $
g_{\alpha}\coloneqq {\rm e}^{\alpha}-{\rm e}^{-\alpha}$.

\subsection{1-point functions}

Also similarly as in \cite{Mason20231-PointAlgebras}, we denote the zero mode of $v\in V_{k}$ with $o(v)\coloneqq v(k-1)$ and extend the definition linearly to every state $u\in V$. Define the 1-point function for $u\in V$ by
\[Z_V(u,\tau)\coloneqq \operatorname{Tr}_Vo(u)q^{L(0)-c/24}.\]
Further, for a module or a twisted module $M$ of a vertex operator algebra $V$, corresponding to a homogeneous state $u\in V_k$, the trace function for $u\in V$ over $M$ can be defined as
\[Z_M(u,\tau)\coloneqq \operatorname{Tr}_Mo(u)q^{L(0)-c/24}.\]

\subsection{Square bracket formalism}

Suppose that $V=\bigoplus_k V_k$ is a vertex operator algebra of central charge $c$. We are going to define some new endomorphisms $v[n]$ for states $v\in V$, called the square bracket modes. The round and square bracket modes are related as follows: for $m\geq 0$ and $v\in V_k$, we have
\begin{gather*}
Y[v,z]=\sum_{n\in \mathbb{Z}}v[n]z^{-n-1}\coloneqq Y\bigl({\rm e}^{zL(0)}v,{\rm e}^z-1\bigr),\\
v[m]=m!\sum_{i\geq m}c(k,i,m)v(i), \qquad m\geq 0,
\end{gather*}
where
\begin{equation}\label{square bracket relations}
 \sum_{m=0}^i c(k,i,m)x^m \coloneqq {k-1+x\choose i}, \qquad \text{where} \ i\geq 0.
\end{equation}

\section[Z\_2-twisted Zhu theory]{$\boldsymbol{Z_2}$-twisted Zhu theory}\label{sec3}

\subsection{Reduction theorems}\label{section3.1}

For the twisted vertex operator defined as in \cite{Frenkel1988VertexMonster},
$
Y \colon V_L\to \operatorname{End}\bigl(V_L^T\bigr)\big[\big[z^{1/2},z^{-1/2}\big]\big]$,
 $z_1,z_2\in \mathbb{C}$, $q_i={\rm e}^{z_i}$, $q={\rm e}^{2\pi{\rm i}\tau}$, define the twisted 2-point function on the torus as
\[F_{V_L^T}((u,z_1),(v,z_2),\tau)\coloneqq \operatorname{Tr}_{V_L^T}Y\bigl(q_1^{L(0)}u,q_1\bigr)Y\bigl(q_2^{L(0)}v,q_2\bigr)q^{L(0)-c/24}.\]
In this subsection, we obtain reduction theorems where we express the above twisted 2-point function on the torus in terms of modular data and 1-point functions.

\begin{Theorem}\label{Reduction Theorem}
For $u,v\in V_L$,
\[
F_{V_L^T}((u,z_1),(v,z_2),\tau)=\sum_{p=0,1}\sum_{m\geq 0}\frac{2^{-m-1}}{m!}Z_{V_L^T}((\theta^pu)[m]v, \tau)P_1^{(m)}(z_{21}/2+p\pi {\rm i},\tau/2),
 \]
 where $z_{21}\coloneqq z_2-z_1=-z_{12}$.
\end{Theorem}

\begin{proof}
Since $L$ is even, for $m,n\in \frac{1}{2}\mathbb{Z}$, we have the following twisted commutator formula from~\cite{Frenkel1988VertexMonster}:
\[[u(m),v(n)]=\frac{1}{2}\sum_{p=0,1}(-1)^{2pm}\sum_{i\in \mathbb{Z}^+}\binom{m}{i}((\theta^pu)(i)v)(m+n-i),\]
and the above equation could further be rewritten as
\[[u(m),Y(v,z)]=\frac{1}{2}\sum_{p=0,1}(-1)^{2pm}\sum_{i\in \mathbb{Z}^+}\binom{m}{i}Y((\theta^pu)(i)v,z)z^{m-i}.\]
For $n\in \frac{1}{2}\mathbb{Z}$, $k\in \mathbb{Z}$, $u\in V_{L_{(k)}}$ (weight $k$ $L(0)$-eigenspace of $V_L$), $v\in V_L$,
\begin{equation}\label{Commutator Formula}
 \big[u(n),Y\bigl(q_2^{L(0)}v, q_2\bigr)\big]=\frac{1}{2}\sum_{p=0,1}(-1)^{2pn}\sum_{i\in \mathbb{Z}^+}\binom{n}{i}Y\bigl((\theta^pu)(i)q_2^{L(0)}v,q_2\bigr)q_2^{n-i},
\end{equation}
and we have
\[\sum_{i\in \mathbb{Z}^+}\binom{n}{i}Y\bigl((\theta^pu)(i)q_2^{L(0)}v,q_2\bigr)q_2^{n-i}=q_2^{n-k+1}Y\biggl(q_2^{L(0)}\sum_{i\in \mathbb{Z}^+}\binom{n}{i}(\theta^pu)(i)v,q_2\biggr).\]
From \eqref{square bracket relations}, we have
\[\sum_{i\in \mathbb{Z}^+}\binom{n}{i}Y\bigl((\theta^pu)(i)q_2^{L(0)}v,q_2\bigr)q_2^{n-i}=q_2^{n-k+1}\sum_{m\geq 0}\frac{(n-k+1)^m}{m!}Y\bigl(q_2^{L(0)} (\theta^p(u))[m]v, q_2\bigr).\]
We set $r\coloneqq n-k+1$, we have
\begin{gather}
\operatorname{Tr}_{V_L^T} \bigl\{ u(n)Y\bigl(q_2^{L(0)}v, q_2\bigr)\bigr\} q^{L(0)-c/24}
\nonumber\\
\qquad= \frac{1}{2} \sum_{p=0,1} (-1)^{2pn}
 \operatorname{Tr}_{V_L^T} q_2^r
 \sum_{m \ge 0} \frac{r^m}{m!}
 Y\bigl(q_2^{L(0)}(\theta^p(u))[m]v, q_2\bigr)
 q^{L(0)-c/24} \nonumber \\
\phantom{\qquad=}{} + \operatorname{Tr}_{V_L^T}
 \bigl\{ Y\bigl(q_2^{L(0)}v, q_2\bigr) u(n) \bigr\}
 q^{L(0)-c/24},
\nonumber
\\
\operatorname{Tr}_{V_L^T} \bigl\{ u(n)Y\bigl(q_2^{L(0)}v, q_2\bigr)\bigr\} q^{L(0)-c/24}\nonumber
\\
\qquad= \frac{1}{2} \sum_{p=0,1} (-1)^{2pn}
 \operatorname{Tr}_{V_L^T} q_2^r
 \sum_{m \ge 0} \frac{r^m}{m!}
 Y\bigl(q_2^{L(0)}(\theta^p(u))[m]v, q_2\bigr)
 q^{L(0)-c/24} \nonumber \\
\phantom{\qquad=}{} + q^r \operatorname{Tr}_{V_L^T} \bigl\{ Y\bigl(q_2^{L(0)}v, q_2\bigr) q^{L(0)-c/24} \bigr\} u(n).\label{eq:trace-commutation}
\end{gather}
Since $\operatorname{Tr}(AB)=\operatorname{Tr}(BA)$, the last term above can be rewritten as follows:
\begin{gather*}
\operatorname{Tr}_{V_L^T} \bigl\{ u(n)Y\bigl(q_2^{L(0)}v, q_2\bigr)\bigr\} q^{L(0)-c/24}
\\
\qquad= \frac{1}{2} \sum_{p=0,1} (-1)^{2pn}
 \operatorname{Tr}_{V_L^T} q_2^r
 \sum_{m \ge 0} \frac{r^m}{m!}
 Y\bigl(q_2^{L(0)}(\theta^p(u))[m]v, q_2\bigr)
 q^{L(0)-c/24} \\
\phantom{\qquad=}{} + q^r \operatorname{Tr}_{V_L^T} \bigl\{ u(n) Y\bigl(q_2^{L(0)}v, q_2\bigr) q^{L(0)-c/24} \bigr\},
\\
\operatorname{Tr}_{V_L^T} \bigl\{ u(n)Y\bigl(q_2^{L(0)}v, q_2\bigr)\bigr\} q^{L(0)-c/24}\\
\qquad= \frac{1}{2} q_2^r \sum_{p=0,1} (-1)^{2pn} \sum_{m \geq 0} \frac{r^m}{m!} \operatorname{Tr}_{V_L^T} Y\bigl(q_2^{L(0)}(\theta^p(u))[m]v, q_2\bigr) q^{L(0)-c/24}\\
\phantom{\qquad=}{} + q^r \operatorname{Tr}_{V_L^T} \bigl\{ u(n) Y\bigl(q_2^{L(0)}v, q_2\bigr) q^{L(0)-c/24} \bigr\},
\\
\operatorname{Tr}_{V_L^T} \bigl\{ u(n)Y\bigl(q_2^{L(0)}v, q_2\bigr)\bigr\} q^{L(0)-c/24}\\
\qquad= \frac{1}{2}q_2^r\sum_{p=0,1}(-1)^{2pn}\sum_{m\geq 0}\frac{r^m}{m!}Z_{V_L^T}((\theta^p(u))[m]v, \tau) \\
\phantom{\qquad=}{}+ q^r\operatorname{Tr}_{V_L^T}\bigl\{u(n)Y\bigl(q_2^{L(0)}v,q_2\bigr)q^{L(0)-c/24}\bigr\}.
\end{gather*}
For $r=0$, the above equation says that
\[
 \sum_{p=0,1}(-1)^{2pn}Z_{V_L^T}((\theta^p(u))[0]v, \tau)=0.
\]
For $r\neq 0$, the above equation gives
\begin{gather}
 \operatorname{Tr}_{V_L^T}\bigl\{ u(n)Y(q_2^{L(0)}v, q_2)\bigr\}q^{L(0)-c/24}\nonumber\\
 \qquad=\frac{q_2^r}{2(1-q^r)}\sum_{p=0,1}(-1)^{2pn}\sum_{m\geq 0}\frac{r^m}{m!}Z_{V_L^T}((\theta^p(u))[m]v, \tau).\label{non-zero r}
\end{gather}
However, observe that the commutator formula \eqref{Commutator Formula} is non-zero only when $n\notin \mathbb{Z}$ (i.e., $r\neq 0$) and hence the 2-point function we defined above can now be written as
\begin{align*}
F_{V_L^T}((u,z_1),(v,z_2),\tau)&= \sum_{n\in \frac{1}{2}\mathbb{Z}}q_1^{k-n-1}\operatorname{Tr}_{V_L^T}\bigl\{u(n)Y\bigl(q_2^{L(0)}v,q_2\bigr)\bigr\}q^{L(0)-c/24}
\\
&= \sum_{n \in \frac{1}{2}\mathbb{Z}, r \neq 0} q_1^{-r} \frac{q_2^r}{2(1-q^r)} \sum_{m \geq 0} \sum_{p=0,1} (-1)^{2pn} \frac{r^m}{m!} Z_{V_L^T}((\theta^p u)[m]v, \tau)
\\
&= \frac{1}{2} \sum_{m \geq 0} \sum_{n \in \frac{1}{2} \mathbb{Z}, r \neq 0} \sum_{p=0,1} \frac{(-1)^{2pn}}{m!} Z_{V_L^T} ((\theta^p u)[m]v, \tau) r^m \frac{q_2^r q_1^{-r}}{1-q^r}.
\end{align*}
We set $q_{21}\coloneqq\frac{q_2}{q_1}$, then
\begin{equation}\label{equation}
 =\frac{1}{2}\sum_{m\geq 0}\sum_{r\in\frac{1}{2}\mathbb{Z},r\neq 0}\sum_{p=0,1}\frac{(-1)^{2pr}}{m!}Z_{V_L^T}((\theta^pu)[m]v, \tau)r^m\frac{q_{21}^r}{(1-q^r)}.
\end{equation}
The above expression can be written as
\[\sum_{p=0,1}\sum_{r\in\frac{1}{2}\mathbb{Z},r\neq 0}(-1)^{2pr}r^m\frac{q_{21}^r}{(1-q^r)}=2^{-m}\sum_{p=0,1}P_1^{(m)}(z_{21}/2+p\pi {\rm i},\tau/2).\]
Thus expression in \eqref{equation} can be rewritten as
\[F_{V_L^T}((u,z_1),(v,z_2),\tau)=\sum_{p=0,1}\sum_{m\geq 0}\frac{2^{-m-1}}{m!}Z_{V_L^T}((\theta^pu)[m]v, \tau)P_1^{(m)}(z_{21}/2+p\pi {\rm i},\tau/2).\tag*{\qed}
\] \renewcommand{\qed}{}
\end{proof}

\begin{Theorem}\label{Reduction Theorem 2}
 The twisted $2$-point function can be rewritten as
 \[\operatorname{Tr}_{V_L^T}Y\bigl(q_1^{L(0)}u,q_1\bigr)Y\bigl(q_2^{L(0)}v,q_2\bigr)q^{L(0)-c/24}=Z_{V_L^T}(Y_{\mathbb{Z}}[u,z_{12}]v,\tau),\]
 where \smash{$Y_{\mathbb{Z}}[v,z]\coloneqq \sum_{n\in \mathbb{Z}}v[n]z^{-n-1}$}, i.e., the restriction of the twisted vertex operator $Y$ {\rm(}which is also referred to as $Y_{V_L^T}$ occasionally in the paper$)$ to integer-modes.
\end{Theorem}

\begin{proof}
 From the associativity \cite[formula~(9.3.52)]{Frenkel1988VertexMonster}, we have
 \[Y_{V_L^T}(u,z_1+z_2)Y_{V_L^T}(v,z_2)=Y_{V_L^T}(Y_{V_L}(u,z_1)v,z_2),\]
 using which we have
 \[\operatorname{Tr}_{V_L^T}Y_{V_L^T}(u,z_1+z_2)Y_{V_L^T}(v,z_2)q^{L(0)-c/24}=\operatorname{Tr}_{V_L^T}Y_{V_L^T}(Y_{V_L}(u,z_1)v,z_2)q^{L(0)-c/24}.\]
 Further, in \cite[Theorem 10]{Mason20231-PointAlgebras} we have the 2-point function $F_V((u,z_1),(v,z_2),\tau)$ rewritten as $Z_V(Y[u,z_{12}]v,\tau)$. Using a similar argument, we have
 \begin{gather*}\operatorname{Tr}_{V_L^T}Y\bigl(q_1^{L(0)}u,q_1\bigr)Y\bigl(q_2^{L(0)}v,q_2\bigr)q^{L(0)-c/24}\\
 \qquad
 =\operatorname{Tr}_{V_L^T}Y\bigl(Y_{\mathbb{Z}}\bigl(q_1^{L(0)}u,q_1-q_2\bigr)q_2^{L(0)}v,q_2\bigr)q^{L(0)-c/24},\end{gather*}
 and thus we obtain the required result.
\end{proof}

\subsection{Recursion formula}\label{sec3.2}

In this subsection, we obtain a recursion formula where we express 1-point functions in terms of modular functions and 1-point functions of states of lower weight.

\begin{Theorem}\label{twisted recursion theorem}
 For $u,v\in V_L^-$, $n\geq 1$,
 \[Z_{V_L^T}(u[-n]v,\tau)=\sum_{m\geq 1}\left(\frac{1}{m}\right)\overline{E}_{m+n}(\tau)Z_{V_L^T}( u[m]v, \tau).\]
\end{Theorem}

\begin{proof}

From Theorems \ref{Reduction Theorem} and \ref{Reduction Theorem 2}, we have
\begin{align*}
Z_{V_L^T}(Y_{\mathbb{Z}}[u, z_{12}] v, \tau) &= \sum_{p=0,1} \sum_{m \geq 1} \frac{2^{-m-1}}{m!} Z_{V_L^T} ((\theta^p u)[m] v, \tau) P_1^{(m)} \left( \frac{z_{21}}{2} + p \pi {\rm i}, \frac{\tau}{2} \right).
\end{align*}
Since
\begin{gather}\label{P_1 translation}
 P_1^{(m)}(z/2+\pi {\rm i},\tau/2)=2^{m+1}P_1^{(m)}(z,\tau)-P_1^{(m)}(z/2,\tau/2),
\end{gather}
we have
\begin{gather*}
Z_{V_L^T}(Y_{\mathbb{Z}}[u,z_{12}]v,\tau)\\
\qquad=\sum_{m\geq 1}\frac{2^{-m-1}}{m!}Z_{V_L^T}((u)[m]v, \tau)P_1^{(m)}(z_{21}/2,\tau/2)\\
\phantom{\qquad=}{}+\sum_{m\geq 1}\frac{2^{-m-1}}{m!}Z_{V_L^T}((\theta u)[m]v, \tau)P_1^{(m)}(z_{21}/2+\pi {\rm i},\tau/2)\\
\qquad=\sum_{m\geq 1}\frac{2^{-m-1}}{m!}Z_{V_L^T}((u)[m]v, \tau)P_1^{(m)}(z_{21}/2,\tau/2)\\
\phantom{\qquad=}{}+\sum_{m\geq 1}\frac{2^{-m-1}}{m!}Z_{V_L^T}((\theta u)[m]v, \tau)\bigl\{ 2^{m+1}P_1^{(m)}(z_{21},\tau)-P_1^{(m)}(z_{21}/2,\tau/2) \bigr\}\\
\qquad=\sum_{m\geq 1}\frac{1}{m!} P_1^{(m)}(z_{21},\tau)Z_{V_L^T}((\theta u)[m]v, \tau)\\
\phantom{\qquad=}{}+\sum_{m\geq 1}\frac{2^{-m-1}}{m!}P_1^{(m)}(z_{21}/2,\tau/2)(Z_{V_L^T}((u)[m]v, \tau)-Z_{V_L^T}((\theta u)[m]v, \tau)).
\end{gather*}
For $n\in \mathbb{Z}^+$, comparing coefficients of $z_{12}^{n-1}$ above gives us
\begin{align*}
Z_{V_L^T}(u[-n]v,\tau)={}&\sum_{m\geq 1}(-1)^{m+1}\binom{m+n-1}{m}E_{m+n}(\tau)Z_{V_L^T}((\theta u)[m]v, \tau)\\
&+\sum_{m\geq 1}(-1)^{m+1}2^{-(m+n)}\binom{m+n-1}{m}E_{m+n}(\tau/2)\\
&\phantom{+}{}\times(Z_{V_L^T}(u[m]v, \tau)-Z_{V_L^T}((\theta u)[m]v, \tau)).\tag*{\qed}
\end{align*}\renewcommand{\qed}{}
\end{proof}

\subsection{Twisted reduction theorems}

For $u,v\in V_L^-$, we can define modified twisted 2-point functions on the torus
\[F_{(V_L^T)^+}((u,z_1),(v,z_2),\tau)\coloneqq \operatorname{Tr}_{(V_L^T)^+}Y\bigl(q_1^{L(0)}u,q_1\bigr)Y\bigl(q_2^{L(0)}v,q_2\bigr)q^{L(0)-c/24}.\]

In this subsection, we obtain expressions for the above modified twisted 2-point functions on the torus in the form of theorems analogous to the reduction theorem in Section~\ref{section3.1}. Before we prove the main theorems, we shall recall the following lemma from \cite[Lemma~12]{Mason20231-PointAlgebras}.

\begin{Lemma}\label{Lemma 12 from MM}
 Let the finite-dimensional linear space $X=X_1\oplus X_2$ decompose as indicated, and let $f,g\in \operatorname{End}(X)$ be a pair of endomorphisms mapping $X_1\to X_2$ and $X_2\to X_1$. Then the following statements hold:
 \begin{enumerate}\itemsep=0pt
 \item[$(1)$] We have $\operatorname{Tr}_{X_2}fg=\operatorname{Tr}_{X_1}gf$.
 \item[$(2)$] We have $\operatorname{Tr}_{X_1}(fg+gf)=\operatorname{Tr}_{X}fg$.
 \item[$(3)$] If $\operatorname{Tr}_X fg=0$ then we have $\operatorname{Tr}_{X_1}fg=-\operatorname{Tr}_{X_1}gf$.
 \end{enumerate}
\end{Lemma}

\begin{Theorem}\label{Z2 Twisted+ Reduction Theorem 1}
 For $u,v\in V_L^-$,
 \begin{gather*}
 F_{(V_L^T)^+}((u,z_1),(v,z_2),\tau)\\
 \qquad=\frac{1}{2}\sum_{p=0,1}\sum_{m\geq 0}\frac{2^{-m}}{m!} Q_1^{(m)}(z_{21}/2+p\pi {\rm i},\tau/2) Z_{(V_L^T)^+}((\theta^pu)[m]v,\tau)\\
 \phantom{\qquad=}{}+\frac{1}{2}\sum_{p=0,1}\sum_{m\geq 1}\frac{2^{-m}}{2(m!)}\bigl( P_1^{(m)}(z_{21}/2+p\pi {\rm i},\tau/2)-Q_1^{(m)}(z_{21}/2+p\pi {\rm i},\tau/2)\bigr) \\
 \phantom{\qquad=+}{}\times Z_{V_L^T}((\theta^p(u))[m]v, \tau).
 \end{gather*}
\end{Theorem}

\begin{proof}

For $u,v\in V_L^-$ and $m,n\in \frac{1}{2}\mathbb{Z}$,
\[
u(m)\colon \ \bigl(V_L^T\bigr)^-\to \bigl(V_L^T\bigr)^+, \qquad v(n)\colon \ \bigl(V_L^T\bigr)^+\to \bigl(V_L^T\bigr)^-,\]
and when $u\in V_{L_{(k)}}$,
\[F_{(V_L^T)^+}((u,z_1),(v,z_2),\tau)\coloneqq \sum_{n\in \frac{1}{2}\mathbb{Z}}q_1^{k-n-1}\operatorname{Tr}_{(V_L^T)^+}\bigl\{u(n)Y\bigl(q_2^{L(0)}v,q_2\bigr)\bigr\}q^{L(0)-c/24}.\]
Using the commutator formula \eqref{Commutator Formula}, and a similar argument as before, where $r\coloneqq n-k+1$, we have
\begin{gather*}
\operatorname{Tr}_{(V_L^T)^+}\bigl\{u(n) Y\bigl(q_2^{L(0)}v,q_2\bigr)\bigr\}q^{L(0)-c/24}
\\
\qquad=\frac{1}{2}q_2^r\sum_{p=0,1}\sum_{m\geq 0}(-1)^{2pr}\frac{r^m}{m!}Z_{(V_L^T)^+}((\theta^p(u))[m]v, \tau)\\
\phantom{\qquad=}{}+\operatorname{Tr}_{(V_L^T)^+}\bigl\{ Y\bigl(q_2^{L(0)}v,q_2\bigr)u(n)\bigr\}q^{L(0)-c/24}\\
\qquad=\frac{1}{2}q_2^r\sum_{p=0,1}\sum_{m\geq 0}(-1)^{2pr}\frac{r^m}{m!}Z_{(V_L^T)^+}((\theta^p(u))[m]v, \tau)\\
\phantom{\qquad=}{}+q^{r}\operatorname{Tr}_{(V_L^T)^+}\bigl\{ Y\bigl(q_2^{L(0)}v,q_2\bigr)q^{L(0)-c/24}u(n)\bigr\}.
\end{gather*}
Using Lemma \ref{Lemma 12 from MM}, we have
\begin{gather*}
=\frac{1}{2}q_2^r\sum_{p=0,1}\sum_{m\geq 0}(-1)^{2pr}\frac{r^m}{m!}Z_{(V_L^T)^+}((\theta^p(u))[m]v, \tau)\\
\qquad{}-q^{r}\operatorname{Tr}_{(V_L^T)^+}\bigl\{u(n) Y\bigl(q_2^{L(0)}v,q_2\bigr)q^{L(0)-c/24}\bigr\}\\
\qquad{}+q^r\operatorname{Tr}_{(V_L^T)}\bigl\{ Y\bigl(q_2^{L(0)}v,q_2\bigr)q^{L(0)-c/24}u(n)\bigr\}.
\end{gather*}
Thus, we have
\begin{gather*}
\operatorname{Tr}_{(V_L^T)^+}\bigl\{u(n) Y\bigl(q_2^{L(0)}v,q_2\bigr)q^{L(0)-c/24}\bigr\}\\
\qquad= \frac{1}{2} \sum_{p=0,1} \frac{q_2^r}{1 + q^r} \sum_{m \geq 0} (-1)^{2pr} \frac{r^m}{m!} Z_{(V_L^T)^+} ( (\theta^p u)[m] v, \tau ) \\
\phantom{\qquad=}{} + \frac{q^r}{1 + q^r} \operatorname{Tr}_{(V_L^T)} \bigl\{ Y\bigl(q_2^{L(0)} v, q_2\bigr) q^{L(0) - c/24} u(n) \bigr\}.
\end{gather*}
Using \eqref{eq:trace-commutation},
\begin{gather*}
=\frac{1}{2}\sum_{p=0,1}\frac{q_2^r}{1+q^r}\sum_{m\geq 0}(-1)^{2pr}\frac{r^m}{m!}Z_{(V_L^T)^+}((\theta^p(u))[m]v, \tau)\\
\qquad{}+\frac{1}{1+q^r}\biggl\{\operatorname{Tr}_{V_L^T}\bigl( u(n)Y\bigl(q_2^{L(0)}v, q_2\bigr)q^{L(0)-c/24}\bigr)\\
\qquad{}-\frac{1}{2}q_2^r\sum_{p=0,1}\sum_{m\geq 0}(-1)^{2pr}\frac{r^m}{m!}Z_{V_L^T}((\theta^p(u))[m]v, \tau)\biggr\}.
\end{gather*}
Thus, we have
\begin{gather*}
\operatorname{Tr}_{(V_L^T)^+}\bigl\{u(n) Y\bigl(q_2^{L(0)}v,q_2\bigr)q^{L(0)-c/24}\bigr\}\\
\qquad=\frac{1}{2}\sum_{p=0,1}\frac{q_2^r}{1+q^r}\sum_{m\geq 0}(-1)^{2pr}\frac{r^m}{m!}\bigl\{Z_{(V_L^T)^+}((\theta^p(u))[m]v, \tau)-Z_{V_L^T}((\theta^p(u))[m]v, \tau)\bigr\}\\
\phantom{\qquad=}{}+\frac{1}{1+q^r}\operatorname{Tr}_{V_L^T}\bigl( u(n)Y\bigl(q_2^{L(0)}v, q_2\bigr)q^{L(0)-c/24}\bigr).
\end{gather*}
Note that when $u,v\in V_L^-$, \smash{$\operatorname{Tr}_{V_L^T}o(u)o(v)q^{L(0)-c/24}=0$} since both $o(u)$ and $o(v)$ are zero as we~set for $\alpha\in \mathfrak{h}$, $\alpha(n)=0$ for $n\in \mathbb{Z}$ (see \cite[equation (9.1.13)]{Frenkel1988VertexMonster}).
Since $u,v\in V_L^-$, we have
\begin{equation}\label{r=0 twisted +}
 \operatorname{Tr}_{(V_L^T)^+}o(u)o(v)q^{L(0)-c/24}=\operatorname{Tr}_{V_L^T}o(u)o(v)q^{L(0)-c/24}=0.
\end{equation}
When $r\neq 0$, using \eqref{non-zero r} we have
\begin{gather*}
\operatorname{Tr}_{(V_L^T)^+}\bigl\{u(n) Y\bigl(q_2^{L(0)}v,q_2\bigr)q^{L(0)-c/24}\bigr\}\\
\qquad
=\frac{1}{2}\sum_{p=0,1}\frac{q_2^r}{1+q^r}\left(\frac{1-q^r}{1-q^r}\right)\sum_{m\geq 0}(-1)^{2pr}\frac{r^m}{m!}\phantom{\qquad
=}{}\\
\phantom{\qquad
=}{}\times\bigl\{Z_{(V_L^T)^+}((\theta^p(u))[m]v, \tau)-Z_{V_L^T}((\theta^p(u))[m]v, \tau)\bigr\}\\
\phantom{\qquad
=\times}{}+\frac{1}{1+q^r}\left(\frac{q_2^r}{2(1-q^r)}\right)\sum_{p=0,1}\sum_{m\geq 1}(-1)^{2pr}\frac{r^m}{m!}Z_{V_L^T}((\theta^p(u))[m]v, \tau).
\end{gather*}
That is, we have
\begin{gather*}
\operatorname{Tr}_{(V_L^T)^+}\bigl\{u(n) Y\bigl(q_2^{L(0)}v,q_2\bigr)q^{L(0)-c/24}\bigr\}\\
\qquad= \frac{1}{2} \sum_{p=0,1} \frac{q_2^r}{1 + q^r} \sum_{m \geq 0} (-1)^{2pr} \frac{r^m}{m!} Z_{(V_L^T)^+} ( (\theta^p u)[m] v, \tau ) \\
\phantom{\qquad
=}{} + \frac{1}{2} \sum_{p=0,1} \frac{q_2^r q^r}{1 - q^{2r}} \sum_{m \geq 1} (-1)^{2pr} \frac{r^m}{m!} Z_{V_L^T} ( (\theta^p u)[m] v, \tau ).
\end{gather*}
Returning to the 2-point function we defined above,
\begin{gather*}
F_{(V_L^T)^+}((u,z_1),(v,z_2),\tau)\\
\qquad= \sum_{n\in \frac{1}{2}\mathbb{Z}}q_1^{k-n-1}\operatorname{Tr}_{(V_L^T)^+}\bigl\{u(n)Y\bigl(q_2^{L(0)}v,q_2\bigr)\bigr\}q^{L(0)-c/24}\\
\qquad=\frac{1}{2}\sum_{p=0,1}\sum_{r\neq 0,r\in \frac{1}{2}\mathbb{Z}}(-1)^{2pr}\frac{q_{21}^r}{1+q^r}\sum_{m\geq 0}\frac{r^m}{m!}Z_{(V_L^T)^+}((\theta^pu)[m]v,\tau)\\
\phantom{\qquad
=}{}+\frac{1}{2}\sum_{p=0,1}\sum_{r\neq 0,r\in \frac{1}{2}\mathbb{Z}}(-1)^{2pr}\frac{q_{21}^rq^r}{1-q^{2r}}\sum_{m\geq 1}\frac{r^m}{m!}Z_{V_L^T}((\theta^p(u))[m]v, \tau)\\
\phantom{\qquad
=}{}+\operatorname{Tr}_{(V_L^T)^+}o(u)o(v)q^{L(0)-c/24}.
\end{gather*}
Using \eqref{r=0 twisted +}, we have
\begin{gather*}
=\frac{1}{2}\sum_{p=0,1}\sum_{r\in \frac{1}{2}\mathbb{Z}}(-1)^{2pr}\frac{q_{21}^r}{1+q^r}\sum_{m\geq 0}\frac{r^m}{m!}Z_{(V_L^T)^+}((\theta^pu)[m]v,\tau)-\frac{1}{2}\sum_{p=0,1}\frac{1}{2}Z_{(V_L^T)^+}((\theta^pu)[0]v,\tau)\\
\phantom{=}{}+\frac{1}{2}\sum_{p=0,1}\sum_{r\neq 0,r\in \frac{1}{2}\mathbb{Z}}(-1)^{2pr}\frac{q_{21}^rq^r}{1-q^{2r}}\sum_{m\geq 1}\frac{r^m}{m!}Z_{V_L^T}((\theta^p(u))[m]v, \tau)\\
=\frac{1}{2}\sum_{p=0,1}\sum_{m\geq 0}\frac{2^{-m}}{m!} Q_1^{(m)}(z_{21}/2+p\pi {\rm i},\tau/2) Z_{(V_L^T)^+}((\theta^pu)[m]v,\tau)\\
\phantom{=}{}-\frac{1}{2}\sum_{p=0,1}\frac{1}{2}Z_{(V_L^T)^+}((\theta^pu)[0]v,\tau)\\
\phantom{=}{}
+\frac{1}{2}\sum_{p=0,1}\frac{1}{2}\sum_{r\neq 0,r\in \frac{1}{2}\mathbb{Z}}(-1)^{2pr}r^mq_{21}^r\left(\frac{1}{1-q^r}-\frac{1}{1+q^r}\right)\sum_{m\geq 1}\frac{1}{m!}Z_{V_L^T}((\theta^p(u))[m]v, \tau)\\
=\frac{1}{2}\sum_{p=0,1}\sum_{m\geq 0}\frac{2^{-m}}{m!} Q_1^{(m)}(z_{21}/2+p\pi {\rm i},\tau/2) Z_{(V_L^T)^+}((\theta^pu)[m]v,\tau)\\
\phantom{=}{}+\frac{1}{2}\sum_{p=0,1}\sum_{m\geq 1}\frac{2^{-m}}{2(m!)}\bigl( P_1^{(m)}(z_{21}/2+p\pi {\rm i},\tau/2)-Q_1^{(m)}(z_{21}/2+p\pi {\rm i},\tau/2)\bigr)\\
\phantom{=+}{}\times Z_{V_L^T}((\theta^p(u))[m]v, \tau).
\end{gather*}
Now, using \eqref{r=0 twisted +}, we have
\begin{gather*}
=\frac{1}{2}\sum_{p=0,1}\sum_{m\geq 0}\frac{2^{-m}}{m!} Q_1^{(m)}(z_{21}/2+p\pi {\rm i},\tau/2) Z_{(V_L^T)^+}((\theta^pu)[m]v,\tau)\\
\phantom{=}{}+\frac{1}{2}\sum_{p=0,1}\sum_{m\geq 1}\frac{2^{-m}}{2(m!)}\bigl( P_1^{(m)}(z_{21}/2+p\pi {\rm i},\tau/2)-Q_1^{(m)}(z_{21}/2+p\pi {\rm i},\tau/2)\bigr)\\
\phantom{=+}{}\times Z_{V_L^T}((\theta^p(u))[m]v, \tau).\tag*{\qed}
\end{gather*}\renewcommand{\qed}{}
\end{proof}

\begin{Theorem}\label{Z2 Twisted+- Reduction Theorem 2}
 For $u,v\in V_L^-$, the modified twisted $2$-point functions can be rewritten as
 \[\operatorname{Tr}_{(V_L^T)^+}Y\bigl(q_1^{L(0)}u,q_1\bigr)Y\bigl(q_2^{L(0)}v,q_2\bigr)q^{L(0)-c/24}=Z_{(V_L^T)^+}(Y_{\mathbb{Z}}[u,z_{12}]v,\tau).\]
\end{Theorem}

\begin{proof}
 Similar to the proof of Theorem \ref{Reduction Theorem 2}.
\end{proof}

\subsection{Twisted recursion formulas}

In this subsection, we obtain recursion formulas for trace functions of states in \smash{$V_L^+$} traced over~\smash{$\bigl(V_L^T\bigr)^+$} analogous to the recursion formula obtained in Section~\ref{sec3.2}. Before we prove the main theorem, we shall recall the following lemma from \cite[Lemma 3]{Mason20231-PointAlgebras}.

\begin{Lemma}\label{Lemma 3 from MM}
 The following hold:
\begin{enumerate}\itemsep=0pt
\item[$(1)$] Laurent series expansions for $P_1(z,\tau)$ and $Q_1(z,\tau)$ are as in \eqref{eq 7 from MM} and \eqref{eq 10 from MM}, respectively.
\item[$(2)$] With respect to the variable $z$, \smash{$P^{(m)}_1(z,\tau)$} and \smash{$Q^{(m)}_1(z,\tau)$} for $m\geq 0$ are odd functions if $m$ is even and even functions if $m$ is odd.
\item[$(3)$] We have $Q_1(z,\tau)=2P_1(z,2\tau)-P_1(z,\tau)$ and $Q_1(z,\tau)$ is an elliptic function for the lattice~${2\pi {\rm i}(\mathbb{Z}\oplus 2\tau\mathbb{Z})}$.
 \end{enumerate}
\end{Lemma}

\begin{Theorem}\label{+-twisted recursion formulas}
 For $u,v\in V_L^-$, $n\geq 1$,
 \begin{align*}
 Z_{(V_L^T)^+}(u[-n]v,\tau)={}&\overline{F}_{0+n}(\tau)Z_{(V_L^T)^+}(u[0]v,\tau)+\sum_{m\geq 1}\frac{1}{m}\overline{F}_{m+n}(\tau)Z_{(V_L^T)^+}(u[m]v,\tau)\\
 &+\frac{1}{2}\sum_{m\geq 1}\frac{1}{m}\left(\overline{E}_{m+n}(\tau)-\overline{F}_{m+n}(\tau)\right)Z_{V_L^T}(u[m]v, \tau).\end{align*}
\end{Theorem}

\begin{proof}

We have
\begin{equation}\label{Q_1 translation}
 Q_1^{(m)}(z/2+\pi {\rm i},\tau/2)=2^{m+1}Q_1^{(m)}(z,\tau)-Q_1^{(m)}(z/2,\tau/2).
\end{equation}
Using Theorem \ref{Z2 Twisted+- Reduction Theorem 2}, we have
\[F_{(V_L^T)^+}((u,z_1),(v,z_2),\tau)=Z_{(V_L^T)^+}(Y_{\mathbb{Z}}[u,z_{12}]v,\tau)\]
and from Theorem \ref{Z2 Twisted+ Reduction Theorem 1},
\begin{gather*}
Z_{(V_L^T)^+}(Y_{\mathbb{Z}}[u,z_{12}]v,\tau)\\
\qquad=\frac{1}{2}\sum_{p=0,1}\sum_{m\geq 0}\frac{2^{-m}}{m!} Q_1^{(m)}(z_{21}/2+p\pi {\rm i},\tau/2) Z_{(V_L^T)^+}((\theta^pu)[m]v,\tau)\\
\phantom{\qquad=}{}
 +\frac{1}{2}\sum_{p=0,1}\sum_{m\geq 1}\frac{2^{-m}}{2(m!)}\bigl( P_1^{(m)}(z_{21}/2+p\pi {\rm i},\tau/2)-Q_1^{(m)}(z_{21}/2+p\pi {\rm i},\tau/2)\bigr)\\
\phantom{\qquad=+}{}\times Z_{V_L^T}((\theta^p(u))[m]v, \tau).
\end{gather*}
Using equations \eqref{P_1 translation}, \eqref{Q_1 translation}, we shall rewrite the above equation as
\begin{gather*}
Z_{(V_L^T)^+}(Y_{\mathbb{Z}}[u,z_{12}]v,\tau)\\
\qquad=\frac{1}{2}\sum_{m\geq 0}\frac{2^{-m}}{m!} Q_1^{(m)}(z_{21}/2,\tau/2) Z_{(V_L^T)^+}(u[m]v,\tau)\\
\phantom{\qquad=}{}+\frac{1}{2}\sum_{m\geq 0}\frac{2^{-m}}{m!} Q_1^{(m)}(z_{21}/2+\pi {\rm i},\tau/2) Z_{(V_L^T)^+}((\theta u)[m]v,\tau)\\
\phantom{\qquad=}{}+\frac{1}{2}\sum_{m\geq 1}\frac{2^{-m}}{2(m!)}\bigl( P_1^{(m)}(z_{21}/2,\tau/2)-Q_1^{(m)}(z_{21}/2,\tau/2)\bigr)Z_{V_L^T}(u[m]v, \tau)\\
\phantom{\qquad=}{}+\frac{1}{2}\sum_{m\geq 1}\frac{2^{-m}}{2(m!)}\bigl( P_1^{(m)}(z_{21}/2+\pi {\rm i},\tau/2)-Q_1^{(m)}(z_{21}/2+\pi {\rm i},\tau/2)\bigr)Z_{V_L^T}((\theta u)[m]v, \tau).\end{gather*}
Thus, we have
\begin{gather*}
Z_{(V_L^T)^+}(Y_{\mathbb{Z}}[u,z_{12}]v,\tau)\\
\qquad=\frac{1}{2}\sum_{m\geq 0}\frac{2^{-m}}{m!} Q_1^{(m)}(z_{21}/2,\tau/2)\{ Z_{(V_L^T)^+}(u[m]v,\tau) - Z_{(V_L^T)^+}((\theta u)[m]v,\tau)\}\\
\phantom{\qquad=}{}+\frac{1}{2}\sum_{m\geq 0}\frac{2}{m!} Q_1^{(m)}(z_{21},\tau) Z_{(V_L^T)^+}((\theta u)[m]v,\tau)\\
\phantom{\qquad=}{}+\frac{1}{2}\sum_{m\geq 1}\frac{2^{-m}}{2(m!)}\bigl( P_1^{(m)}(z_{21}/2,\tau/2)-Q_1^{(m)}(z_{21}/2,\tau/2)\bigr)\\
\phantom{\qquad=+}{}\times\{Z_{V_L^T}(u[m]v, \tau)-Z_{V_L^T}((\theta u)[m]v, \tau)\}\\
\phantom{\qquad=}{}+\frac{1}{2}\sum_{m\geq 1}\frac{2}{2(m!)}\bigl( P_1^{(m)}(z_{21},\tau)-Q_1^{(m)}(z_{21},\tau)\bigr)Z_{V_L^T}((\theta u)[m]v, \tau).\end{gather*}
Now using Lemma \ref{Lemma 3 from MM}, we shall rewrite the above as
\begin{gather*}
Z_{(V_L^T)^+}(Y_{\mathbb{Z}}[u,z_{12}]v,\tau)\\
\qquad=\frac{1}{2}\sum_{m\geq 0}\frac{(-1)^{m+1}2^{-m}}{m!} Q_1^{(m)}(z_{12}/2,\tau/2)\{ Z_{(V_L^T)^+}(u[m]v,\tau) - Z_{(V_L^T)^+}((\theta u)[m]v,\tau)\}\\
\phantom{\qquad=}{}+\frac{1}{2}\sum_{m\geq 0}\frac{(-1)^{m+1}2}{m!} Q_1^{(m)}(z_{12},\tau) Z_{(V_L^T)^+}((\theta u)[m]v,\tau)\\
\phantom{\qquad=}{}+ \frac{1}{2} \sum_{m \geq 1} \frac{(-1)^{m+1} 2^{-m}}{2(m!)} \bigl( P_1^{(m)}(z_{12}/2, \tau/2) - Q_1^{(m)}(z_{12}/2, \tau/2) \bigr)\\
\phantom{\qquad=+}{}\times \{ Z_{V_L^T}(u[m] v, \tau) - Z_{V_L^T} ( (\theta u)[m] v, \tau ) \}
\\
\phantom{\qquad=}{}+\frac{1}{2}\sum_{m\geq 1}\frac{(-1)^{m+1}2}{2(m!)}\bigl( P_1^{(m)}(z_{12},\tau)-Q_1^{(m)}(z_{12},\tau)\bigr)Z_{V_L^T}((\theta u)[m]v, \tau).
\end{gather*}
Comparing the coefficients of $z_{12}^{n-1}$ above (where $n\in \mathbb{Z}^+$), we have
\begin{gather*}
Z_{(V_L^T)^+}(u[-n]v,\tau)\\
\qquad=\frac{1}{2}\sum_{m\geq 0}(-1)^{m+1}2^{-m-n+1}\binom{m+n-1}{m}F_{m+n}(\tau/2)\\
\phantom{\qquad=+}{}\times \{ Z_{(V_L^T)^+}(u[m]v,\tau) - Z_{(V_L^T)^+}((\theta u)[m]v,\tau)\}\\
\phantom{\qquad=}{}+\frac{1}{2}\sum_{m\geq 0}(-1)^{m+1}2\binom{m+n-1}{m}F_{m+n}(\tau) Z_{(V_L^T)^+}((\theta u)[m]v,\tau)\\
\phantom{\qquad=}{}
+ \frac{1}{2} \sum_{m \geq 1} \frac{(-1)^{m+1} 2^{-m-n+1}}{2} \binom{m+n-1}{m} ( E_{m+n}(\tau/2) - F_{m+n}(\tau/2) ) \\
\phantom{\qquad=+}{}\times \{ Z_{V_L^T}((u)[m]v, \tau) - Z_{V_L^T}((\theta u)[m]v, \tau) \}
\\
\phantom{\qquad=}{}+\frac{1}{2}\sum_{m\geq 1}\frac{(-1)^{m+1}2}{2}\binom{m+n-1}{m}(E_{m+n}(\tau)-F_{m+n}(\tau))Z_{V_L^T}((\theta u)[m]v, \tau).
\end{gather*}
Since $u\in V_L^-$, we have $\theta u=-u$, and hence
\begin{gather*}
Z_{(V_L^T)^+}(u[-n]v,\tau)\\
\qquad=\sum_{m\geq 0}(-1)^{m+1}2^{-m-n+1}\binom{m+n-1}{m}F_{m+n}(\tau/2) Z_{(V_L^T)^+}(u[m]v,\tau)\\
\phantom{\qquad=}{}-\sum_{m\geq 0}(-1)^{m+1}\binom{m+n-1}{m}F_{m+n}(\tau) Z_{(V_L^T)^+}(u[m]v,\tau)\\
\phantom{\qquad=}{}+\frac{1}{2}\sum_{m\geq 1}(-1)^{m+1}2^{-m-n+1}\binom{m+n-1}{m}(E_{m+n}(\tau/2)-F_{m+n}(\tau/2))Z_{V_L^T}(u[m]v, \tau)\\
\phantom{\qquad=}{}-\frac{1}{2}\sum_{m\geq 1}(-1)^{m+1}\binom{m+n-1}{m}(E_{m+n}(\tau)-F_{m+n}(\tau))Z_{V_L^T}(u[m]v, \tau).
\end{gather*}
Using the notation for renormalized Eisenstein series, we have
\begin{align*}
Z_{(V_L^T)^+}(u[-n]v,\tau)={}&\sum_{m\geq 1}\frac{1}{m}\bigl\{2^{-m-n+1}\hat{F}_{m+n}(\tau/2) -\hat{F}_{m+n}(\tau)\bigr\}Z_{(V_L^T)^+}(u[m]v,\tau)\\
&+\bigl\{2^{-n+1}\hat{F}_{0+n}(\tau/2) -\hat{F}_{0+n}(\tau)\bigr\}Z_{(V_L^T)^+}(u[0]v,\tau)\\
&+\frac{1}{2}\sum_{m\geq 1}\frac{1}{m}\bigl\{2^{-m-n+1}\hat{E}_{m+n}(\tau/2) -\hat{E}_{m+n}(\tau)\bigr\}Z_{V_L^T}(u[m]v, \tau)\\
&-\frac{1}{2}\sum_{m\geq 1}\frac{1}{m}\bigl\{2^{-m-n+1}\hat{F}_{m+n}(\tau/2) -\hat{F}_{m+n}(\tau)\bigr\}Z_{V_L^T}(u[m]v, \tau).
\end{align*}
Thus, the above can rewritten as
\begin{align*}
Z_{(V_L^T)^+}(u[-n]v,\tau)={}& \overline{F}_{0+n}(\tau)Z_{(V_L^T)^+}(u[0]v,\tau)+\sum_{m\geq 1}\frac{1}{m}\overline{F}_{m+n}(\tau)Z_{(V_L^T)^+}(u[m]v,\tau)\\
& +\frac{1}{2}\sum_{m\geq 1}\frac{1}{m}\bigl(\overline{E}_{m+n}(\tau)-\overline{F}_{m+n}(\tau)\bigr)Z_{V_L^T}(u[m]v, \tau).\tag*{\qed}
\end{align*}\renewcommand{\qed}{}
\end{proof}

\section{Calculations in the twisted space}\label{sec4}

In this section, we perform certain computations as applications of the theorems we obtained in the previous section.

\subsection{In the twisted module}

Here we derive formulas corresponding to the twisted space using the recursion formula we obtained in Section~\ref{sec3.2}.

\begin{Theorem}\label{Calculation in Twisted Module}
 For positive integers $n_j\geq 1$, corresponding to the state
$u\coloneqq h_{i_1}[-n_1]\cdots\allowbreak h_{i_p}[-n_p]e_{\alpha}$,
 where $e_{\alpha}=f_{\alpha}$ or $g_{\alpha}$ according as $p$ is even or odd respectively, we have
 \[Z_{(V_L^T)}(u,\tau)=\bigg\{\sum_{\sigma}\prod_{(rs)(t)}\delta_{i_r,i_s} \Tilde{\overline{E}}_{n_t}\overline{E}_{n_r+n_s}(\tau)\bigg\}Z_{(V_L^T)}(f_{\alpha},\tau),\]
 where
 $\Tilde{\overline{E}}_{n_k}(\tau)\coloneqq \langle h_{i_k},\alpha \rangle \overline{E}_{n_k}(\tau)$
 where $\langle \cdot,\cdot\rangle$ is naturally extended from $L$ to $\mathfrak{h}=L\otimes \mathbb{C}$ and $\sigma$ ranges over all involutions of the index set $S_p=\{1,2,\dots, p\}$ and $(rs)$, $(t)$ range over the $2$-cycles and $1$-cycles respectively in the decomposition of $\sigma$ in $S_p$ as a product of disjoint $2$-cycles and $1$-cycles.
\end{Theorem}

\begin{proof}
 Let
$v\coloneqq h_{i_2}[-n_2]\cdots h_{i_p}[-n_p]e_{\alpha}$,
 then Theorem \ref{twisted recursion theorem} gives
 \[Z_{(V_L^T)}(u,\tau)=\sum_{m\geq 1}\left(\frac{1}{m}\right) \overline{E}_{m+n}(\tau) Z_{V_L^T}( h_{i_1}[m]v, \tau).\]
Since
 \[\sum_{m\geq 1}h_{i_1}[m]v=\sum_{j=2}^p \delta_{i_1,i_j}n_j (v \backslash j),\]
where $v\backslash j$ denotes the state obtained from $v$ by deleting the operator $h_{i_j}[-n_j]$, we have
 \begin{align*}
 Z_{V_L^T}(u,\tau)={}&\sum_{j=2}^p\delta_{i_1,i_j}\overline{E}_{n_j+n_1}(\tau)Z_{V_L^T}( v\backslash j, \tau)\\
 &+\langle h_{i_1},\alpha\rangle \bigl\{2^{-(n_1-1)}{E}_{n_1}(\tau/2)-{E}_{n_1}(\tau)\bigr\} Z_{V_L^T}( v', \tau),\end{align*}
 where $v'$ denotes the state
$v'\coloneqq h_{i_2}[-n_2]\cdots h_{i_p}[-n_p]e_{\alpha}'$
 and $e_{\alpha}'=g_{\alpha}$ or $f_{\alpha}$ according as $e_{\alpha}=f_{\alpha}$ or $g_{\alpha}$, respectively.
Thus similarly as in \cite{Mason20231-PointAlgebras} one can easily obtain using induction that
\[
Z_{(V_L^T)}(u,\tau)=\biggl(\sum_{\sigma}\prod_{(rs)(t)}\delta_{i_r,i_s}
 \Tilde{\overline{E}}_{n_t}\overline{E}_{n_r+n_s}(\tau)\biggr)Z_{V_L^T}(f_{\alpha},\tau).\tag*{\qed}
\]\renewcommand{\qed}{}
\end{proof}

An immediate corollary of the above theorem is the following.

\begin{Corollary}\label{corollary1}
 For positive integers $n_j\geq 1$, and even $p$, corresponding to the state
$
u\coloneqq h_{i_1}[-n_1]\cdots h_{i_p}[-n_p]\mathbf{1}$,
 we have
 \[Z_{(V_L^T)}(u,\tau)=\bigg\{\sum_{\sigma}\prod_{(rs)}\delta_{i_r,i_s}\overline{E}_{n_r+n_s}(\tau)\bigg\}Z_{(V_L^T)}(\mathbf{1},\tau),\]
 where $\sigma$ ranges over all fixed point free involutions of $S_p$ and $(rs)$ ranges over all $2$-cycles in the decomposition of $\sigma$ in $S_p$.
\end{Corollary}

\begin{proof}
 Evaluate the above theorem at $\alpha=0$.
\end{proof}

\subsection{Outside 2L}

Here we obtain formulas corresponding to certain states in $V_L^+$ whose lattice part comes from an element outside $2L$.

\begin{Theorem}\label{Outside 2L}
 If $\alpha\in L\backslash 2L$, for positive integers $n_j\geq 1$, corresponding to the state
 \[u\coloneqq h_{i_1}[-n_1]\cdots h_{i_p}[-n_p]e_{\alpha},\]
 where $e_{\alpha}=f_{\alpha}$ or $g_{\alpha}$ according as $p$ is even or odd, respectively, we have
 \[Z_{(V_L^T)^+}(u,\tau)=\bigg\{\sum_{\sigma} \prod_{(rs)(t)}\delta_{i_r,i_s}\Tilde{\overline{F}}_{n_t}(\tau)\overline{F}_{n_r+n_s}(\tau)\bigg\}Z_{(V_L^T)^+}(f_{\alpha},\tau)\]
 where
\smash{$\Tilde{\overline{F}}_{n_k}(\tau)\coloneqq \langle h_{i_k},\alpha \rangle\overline{F}_{n_k}(\tau)$}
 where $\langle \cdot,\cdot \rangle$ is naturally extended from $L$ to $\mathfrak{h}=L\otimes \mathbb{C}$ and $\sigma$ ranges over all involutions of the index set $S_p=\{1,2,\dots, p\}$ and $(rs)$, $(t)$ range over the $2$-cycles and $1$-cycles respectively in the decomposition of $\sigma$ in $S_p$ as a product of disjoint $2$-cycles and $1$-cycles.
\end{Theorem}

\begin{proof}
Let
$v\coloneqq h_{i_2}[-n_2]\cdots h_{i_p}[-n_p]e_{\alpha}$,
$v'\coloneqq h_{i_2}[-n_2]\cdots h_{i_p}[-n_p]e_{\alpha}'$,
where $e_{\alpha}'\coloneqq g_{\alpha}$ or $f_{\alpha}$ according as $e_{\alpha}=f_{\alpha}$ or $g_{\alpha}$, respectively.
From Theorem \ref{+-twisted recursion formulas}, we have
 \begin{align*}
 Z_{(V_L^T)^+}(u[-n]v,\tau)={}&\overline{F}_{0+n}(\tau)Z_{(V_L^T)^+}(u[0]v,\tau)+\sum_{m\geq 1}\frac{1}{m}\overline{F}_{m+n}(\tau)Z_{(V_L^T)^+}(u[m]v,\tau)\\
 &+\frac{1}{2}\sum_{m\geq 1}\frac{1}{m}\left(\overline{E}_{m+n}(\tau)-\overline{F}_{m+n}(\tau)\right)Z_{V_L^T}(u[m]v, \tau).
 \end{align*}
Thus,
 \begin{align*}
 Z_{(V_L^T)^+}(u,\tau)={}&\overline{F}_{n_1+0}(\tau)Z_{(V_L^T)^+}(h_{i_1}[0]v,\tau)+\sum_{m\geq 1}\left(\frac{1}{m}\right)\overline{F}_{m+n_1}(\tau)Z_{(V_L^T)^+}(h_{i_1}[m]v,\tau)\\
&+\frac{1}{2}\sum_{m\geq 1}\left(\frac{1}{m}\right)\left\{\overline{E}_{m+n_1}(\tau)-\overline{F}_{m+n_1}(\tau)\right\}Z_{V_L^T}(h_{i_1}[m]v, \tau).
\end{align*}
If $\alpha\notin 2L$, since $Z_{V_L^T}(h_{i_1}[m]v, \tau)=0$ due to the structure of $T$, we have
\begin{gather*}
Z_{(V_L^T)^+}(u,\tau)=\overline{F}_{n_1+0}(\tau)Z_{(V_L^T)^+}(h_{i_1}[0]v,\tau)+\sum_{m\geq 1}\left(\frac{1}{m}\right)\overline{F}_{m+n_1}(\tau)Z_{(V_L^T)^+}(h_{i_1}[m]v,\tau),
\\
Z_{(V_L^T)^+}(u,\tau)=\bigl\{2^{-(n_1-1)}{F}_{n_1}(\tau/2) -{F}_{n_1}(\tau)\bigr\}Z_{(V_L^T)^+}(h_{i_1}[0]v,\tau)\\
\phantom{Z_{(V_L^T)^+}(u,\tau)=}{}+\sum_{j=2}^p\left(\frac{1}{n_j}\right)\bigl\{2^{-(n_j+n_1-1)}\hat{F}_{n_j+n_1}(\tau/2) -\hat{F}_{n_j+n_1}(\tau)\bigr\}Z_{(V_L^T)^+}(h_{i_1}[n_j]v,\tau),\\
Z_{(V_L^T)^+}(u,\tau)=\langle h_{i_1},\alpha\rangle\bigl\{2^{-(n_1-1)}{F}_{n_1}(\tau/2) -{F}_{n_1}(\tau)\bigr\}Z_{(V_L^T)^+}(v',\tau)\\
\phantom{Z_{(V_L^T)^+}(u,\tau)=}{}+\sum_{j=2}^p\delta_{i_1,i_j}\bigl\{2^{-(n_1+n_j-1)}\hat{F}_{n_j+n_1}(\tau/2) -\hat{F}_{n_j+n_1}(\tau)\bigr\}Z_{(V_L^T)^+}(v\backslash j,\tau),
\end{gather*}
which can further be rewritten as
\[Z_{(V_L^T)^+}(u,\tau)=\Tilde{\overline{F}}_{n_1}(\tau)Z_{(V_L^T)^+}(v',\tau)+\sum_{j=2}^p \delta_{i_1,i_j}\overline{F}_{n_j+n_1}(\tau)Z_{(V_L^T)^+}(v\backslash j,\tau).\]
Thus, we have
\[
Z_{(V_L^T)^+}(u,\tau)=\bigg\{\sum_{\sigma} \prod_{(rs)(t)}\delta_{i_r,i_s}\Tilde{\overline{F}}_{n_t}(\tau)
\overline{F}_{n_r+n_s}(\tau)\bigg\}Z_{(V_L^T)^+}(f_{\alpha},\tau).\tag*{\qed}
\]\renewcommand{\qed}{}
\end{proof}

\subsection{At Heisenberg states}

Here we obtain formulas corresponding to the Heisenberg states in $V_L^+$.

\begin{Theorem}\label{At Heisenberg states}
 For positive integers $n_j\geq 1$ and $p$ even, corresponding to the state
$u\coloneqq h_{i_1}[-n_1]\cdots h_{i_p}[-n_p]\mathbf{1}$,
 we have
 \[Z_{(V_L^T)^+}(u,\tau)=\frac{1}{2}Z_{V_L^T}(u,\tau)+\sum_{\sigma\in \operatorname{Inv}_0(\underline{p})}\prod_{(rs)}\delta_{i_r,i_s}\overline{F}_{n_r+n_s}(\tau)\left\{Z_{(V_L^T)^+}(\mathbf{1},\tau)-\frac{1}{2}Z_{V_L^T}(\mathbf{1},\tau)\right\},\]
 where $\underline{p}\coloneqq \{1,2,\dots, p\}$ and $\operatorname{Inv}_0(\underline{p})$ is the set of all fixed-point-free involutions of $\underline{p}$.
\end{Theorem}

\begin{proof}

From Theorem \ref{+-twisted recursion formulas}, if
$v\coloneqq h_{i_2}[-n_2]\cdots h_{i_p}[-n_p]\mathbf{1}$,
then for $n=n_1$, we have
\begin{align*}
Z_{(V_L^T)^+}(u,\tau)={}&\sum_{m\geq 1}\frac{1}{m}\overline{F}_{m+n}(\tau)Z_{(V_L^T)^+}(h_{i_1}[m]v,\tau)+\frac{1}{2}\sum_{m\geq 1}\bigl(\overline{E}_{m+n}(\tau)-\overline{F}_{m+n}(\tau)\bigr)\\
&\times Z_{V_L^T}(h_{i_1}[m]v, \tau).
\end{align*}
Using Heisenberg relations, we have
\begin{gather*}
Z_{(V_L^T)^+}(u,\tau)=\sum_{j=2}^p\delta_{i_1,i_j}\overline{F}_{n_1+n_j}(\tau)\left\{Z_{(V_L^T)^+}(v\backslash j,\tau)-\frac{1}{2}Z_{V_L^T}(v\backslash j,\tau)\right\}\\
\phantom{Z_{(V_L^T)^+}(u,\tau)=}{}+\frac{1}{2}\sum_{j=2}^p\delta_{i_1,i_j}\overline{E}_{n_1+n_j}(\tau)Z_{V_L^T}(v\backslash j,\tau),\\
Z_{(V_L^T)^+}(u,\tau)-\frac{1}{2}Z_{V_L^T}(u,\tau)=\sum_{j=2}^p\delta_{i_1,i_j}\overline{F}_{n_1+n_j}(\tau)\left\{Z_{(V_L^T)^+}(v\backslash j,\tau)-\frac{1}{2}Z_{V_L^T}(v\backslash j,\tau)\right\}.
\end{gather*}
Using induction further, we obtain
\[
Z_{(V_L^T)^+}(u,\tau)-\frac{1}{2}Z_{V_L^T}(u,\tau)=\sum_{\sigma\in {\rm Inv}_0(\underline{p})}\prod_{(rs)}\delta_{i_r,i_s}\overline{F}_{n_r+n_s}(\tau)
\left\{Z_{(V_L^T)^+}(\mathbf{1},\tau)-\frac{1}{2}Z_{V_L^T}(\mathbf{1},\tau)\right\}.\tag*{\qed}
\]\renewcommand{\qed}{}
\end{proof}

\subsection{Inside 2L}

Here we obtain formulas corresponding to certain states in $V_L^+$ whose lattice part comes from an element inside $2L$.

\begin{Theorem}\label{Inside 2L}
 For positive integers $n_j\geq 1$, $\alpha\in 2L$, corresponding to the state
$u\coloneqq h_{i_1}[-n_1]\cdots\allowbreak h_{i_p}[-n_p]e_{\alpha}$,
 where $e_{\alpha}$ is $f_{\alpha}$ or $g_{\alpha}$ according as $p$ is even or odd respectively, we have
 \begin{align*}
 Z_{(V_L^T)^+}(u,\tau)={}&\biggl(\sum_{\sigma}\prod_{(rs)(t)}\delta_{i_r,i_s}\Tilde{\overline{F}}_{n_t}(\tau)\overline{F}_{n_r+n_s}(\tau)\biggr)
 \left\{Z_{(V_L^T)^+}(f_{\alpha},\tau)-\frac{1}{2}Z_{V_L^T}(f_{\alpha},\tau)\right\}\\
 &+\frac{1}{2}Z_{V_L^T}(u,\tau),
 \end{align*}
 where $\sigma$ ranges over all involutions of the index set $S_p=\{1,2,\dots, p\}$ and $(rs)$, $(t)$ range over the $2$-cycles and $1$-cycles respectively in the decomposition of $\sigma$ in $S_p$ as a product of disjoint $2$-cycles and $1$-cycles.
\end{Theorem}

\begin{proof}

For $u,v\in V_L^-$, from Theorem \ref{+-twisted recursion formulas}, we have
 \begin{align*}
 Z_{(V_L^T)^+}(u[-n]v,\tau)={}&\overline{F}_{0+n}(\tau)Z_{(V_L^T)^+}(u[0]v,\tau)+\sum_{m\geq 1}\frac{1}{m}\overline{F}_{m+n}(\tau)Z_{(V_L^T)^+}(u[m]v,\tau)\\
 &+\frac{1}{2}\sum_{m\geq 1}\frac{1}{m}\left(\overline{E}_{m+n}(\tau)-\overline{F}_{m+n}(\tau)\right)Z_{V_L^T}(u[m]v, \tau).
 \end{align*}
Let
$v\coloneqq h_{i_2}[-n_2]\cdots h_{i_p}[-n_p]e_{\alpha}$,
$v'\coloneqq h_{i_2}[-n_2]\cdots h_{i_p}[-n_p]e'_{\alpha}$,
where we use same notation as earlier for $e_{\alpha}'$. Then similarly as earlier, we have
 \begin{align*}
 Z_{(V_L^T)^+}(u,\tau)={}&\sum_{j=2}^p\delta_{i_1,i_j}\overline{F}_{n_1+n_j}(\tau)\left\{Z_{(V_L^T)^+}(v\backslash j,\tau)-\frac{1}{2}Z_{V_L^T}(v\backslash j,\tau)\right\}\\
&+\frac{1}{2}\sum_{j=2}^p\delta_{i_1,i_j}\overline{E}_{n_1+n_j}(\tau)Z_{V_L^T}(v\backslash j,\tau)+\Tilde{\overline{F}}_{n_1}(\tau)Z_{(V_L^T)^+}(v',\tau),
\end{align*}
which can be rewritten as
 \begin{align*}
 Z_{(V_L^T)^+}(u,\tau)-\frac{1}{2}Z_{V_L^T}(u,\tau)={}&\sum_{j=2}^p\delta_{i_1,i_j}\overline{F}_{n_1+n_j}(\tau)\left\{Z_{(V_L^T)^+}(v\backslash j,\tau)-\frac{1}{2}Z_{V_L^T}(v\backslash j,\tau)\right\}\\
&+\Tilde{\overline{F}}_{n_1}(\tau)Z_{(V_L^T)^+}(v',\tau).
\end{align*}
Due to the structure of $T$, we have
$Z_{V_L^T}(v',\tau)=0$.
Thus, we can rewrite the above equation as
 \begin{align*}
 Z_{(V_L^T)^+}(u,\tau)-\frac{1}{2}Z_{V_L^T}(u,\tau)={}&\sum_{j=2}^p\delta_{i_1,i_j}\overline{F}_{n_1+n_j}(\tau)\left\{Z_{(V_L^T)^+}(v\backslash j,\tau)-\frac{1}{2}Z_{V_L^T}(v\backslash j,\tau)\right\}\\
&+\Tilde{\overline{F}}_{n_1}(\tau)\left\{Z_{(V_L^T)^+}(v',\tau)-\frac{1}{2}Z_{V_L^T}(v',\tau)\right\}.
\end{align*}
We thus have
\begin{gather*}
Z_{(V_L^T)^+}(u,\tau)-\frac{1}{2}Z_{V_L^T}(u,\tau)\\
\qquad=\biggl(\sum_{\sigma}\prod_{(rs)(t)}\delta_{i_r,i_s}\Tilde{\overline{F}}_{n_t}
(\tau)\overline{F}_{n_r+n_s}(\tau)\biggr)\left\{Z_{(V_L^T)^+}(f_{\alpha},\tau)
-\frac{1}{2}Z_{V_L^T}(f_{\alpha},\tau)\right\}.\tag*{\qed}
\end{gather*}\renewcommand{\qed}{}
\end{proof}

\section{1-point functions}\label{sec5}

In this section, we write down the explicit formulas for the 1-point functions of all states in the $\mathbb{Z}_2$-orbifold of the lattice vertex operator algebra. Recall that $l=k/8$ where $k$ is the rank of the lattice $L$.

\begin{Lemma}\label{vacuum traces}
 We have
 \begin{gather}\label{vacuum trace of lattice VOA}
 Z_{V_L}(\mathbf{1},\tau)=\frac{\theta_L(\tau)}{\eta(\tau)^k},
\\
\label{vacuum trace of twisted module of lattice VOA}
 Z_{V_L^T}(\mathbf{1},\tau)=\eta(\tau)^{k/2}\left(\frac{\Theta_2(\tau)}{2}\right)^{-k/2},
\\
\label{vacuum trace over first summand}
 Z_{V_L^+}(\mathbf{1},\tau)=\frac{1}{2}\frac{\theta_L(\tau)}{\eta(\tau)^k}+\frac{1}{2}\eta(\tau)^{k/2}\left(\frac{\Theta_1(\tau)}{2}\right)^{-k/2},
\\
\label{vacuum trace over second summand}
 Z_{(V_L^T)^+}(\mathbf{1},\tau)=\frac{1}{2}\eta(\tau)^{k/2} \left\{\left(\frac{\Theta_2(\tau)}{2}\right)^{-k/2}+(-1)^l\left(\frac{\Theta_3(\tau)}{2}\right)^{-k/2}\right\}.
 \end{gather}
\end{Lemma}

\begin{proof}
 We know equations \eqref{vacuum trace of lattice VOA}, \eqref{vacuum trace over first summand} from \cite[Lemma 7]{Mason20231-PointAlgebras} and equations \eqref{vacuum trace of twisted module of lattice VOA} and \eqref{vacuum trace over second summand} can be obtained using similar counting.
\end{proof}

\begin{Lemma}\label{lattice element traces}
For $\alpha\in L\backslash 2L$, we have
\smash{$
Z_{V_L^+}(f_{\alpha},\tau)=0$}, \smash{$ Z_{(V_L^T)^+}(f_{\alpha},\tau)=0$}.
For $\alpha\in 2L$,
\begin{gather*}
 Z_{V_L^+}(f_{\alpha},\tau)=\eta(\tau)^{k/2}\left(\frac{\Theta_1(\tau)}{2}\right)^{\langle\alpha,\alpha\rangle-k/2},
\\
 \frac{1}{2}Z_{V_L^T}(f_{\alpha},\tau)=\eta(\tau)^{k/2}\left\{\left(\frac{\Theta_2(\tau)}{2}\right)^{\langle\alpha,\alpha\rangle-k/2}\right\},
\\
 Z_{(V_L^T)^+}(f_{\alpha},\tau)=\eta(\tau)^{k/2}\left\{\left(\frac{\Theta_2(\tau)}{2}\right)^{\langle\alpha,\alpha\rangle-k/2}+(-1)^l\left(\frac{\Theta_3(\tau)}{2}\right)^{\langle\alpha,\alpha\rangle-k/2}\right\}.
\end{gather*}
\end{Lemma}

\begin{proof}
 Similarly as in \cite{Dong2000MonstrousWeight}.
\end{proof}

We have the following formulas for 1-point functions in the $\mathbb{Z}_2$-orbifold $V$ of $V_L$ given by
\[V=V_L^+\oplus \bigl(V_L^T\bigr)^+.\]

\begin{Lemma}
 Suppose $V$ is the $\mathbb{Z}_2$-orbifold of the lattice VOA $V_L$ corresponding to the $(-1)$-involution as defined in Section {\rm\ref{sec2.3}}, then we have
 \[Z_V(\mathbf{1},\tau)=\frac{1}{2}\frac{\theta_L(\tau)}{\eta(\tau)^{k}}+\frac{1}{2}\eta(\tau)^{k/2}\left\{ \left(\frac{\Theta_1(\tau)}{2}\right)^{-k/2}+ \left(\frac{\Theta_2(\tau)}{2}\right)^{-k/2}+(-1)^l \left(\frac{\Theta_3(\tau)}{2}\right)^{-k/2} \right\}\]
 and for $\alpha\in 2L$,
 \[Z_V(f_{\alpha},\tau)=\eta(\tau)^{k/2}\left\{\left(\frac{\Theta_1(\tau)}{2}\right)^{\langle\alpha,\alpha\rangle-k/2}\!
 +\left(\frac{\Theta_2(\tau)}{2}\right)^{\langle\alpha,\alpha\rangle-k/2}\!+(-1)^l \left(\frac{\Theta_3(\tau)}{2}\right)^{\langle\alpha,\alpha\rangle-k/2}\right\}\]
 for $\alpha\in L\backslash 2L$,
$Z_V(f_{\alpha},\tau)=0$.
\end{Lemma}

\begin{proof}
 Using Lemmas \ref{vacuum traces} and~\ref{lattice element traces}, this follows.
\end{proof}

\begin{Theorem}

Suppose $V$ is the $\mathbb{Z}_2$-orbifold of the lattice VOA $V_L$ corresponding to the $(-1)$-involution as defined in Section {\rm\ref{sec2.3}}, then for positive integers $n_i\geq 1$, corresponding to the Heisenberg state
$u\coloneqq h_{i_1}[-n_1]\cdots h_{i_p}[-n_p]\mathbf{1}\in V$,
we have
\begin{align*}
Z_V(u,\tau)={}&\frac{1}{2}\sum_{\Delta\subseteq \Lambda}\frac{\theta_L(\tau,P_{\Delta})}{\eta(\tau)^k}\bigg(\sum_{\sigma\in \operatorname{Inv}_0(\underline{p}\backslash \Delta)}\prod_{(rs)}\delta_{i_r,i_s}\hat{E}_{n_r+n_s}(\tau)\bigg)\\
&+\frac{1}{2}\eta(\tau)^{k/2}\left(\frac{\Theta_1(\tau)}{2}\right)^{-k/2}\bigg(\sum_{\sigma\in \operatorname{Inv}_0(\underline{p})} \prod_{(rs)}\delta_{i_r,i_s}\hat{F}_{n_r+n_s}(\tau) \bigg)\\
&+\frac{1}{2}\eta(\tau)^{k/2}\left(\frac{\Theta_2(\tau)}{2}\right)^{-k/2}\bigg(\sum_{\sigma\in \operatorname{Inv}_0(\underline{p})}\prod_{(rs)}\delta_{i_r,i_s}\overline{E}_{n_r+n_s}(\tau)\bigg)\\
&+\frac{(-1)^l}{2}\eta(\tau)^{k/2}\left(\frac{\Theta_3(\tau)}{2}\right)^{-k/2}\bigg(\sum_{\sigma\in \operatorname{Inv}_0(\underline{p})}\prod_{(rs)}\delta_{i_r,i_s}\overline{F}_{n_r+n_s}(\tau)\bigg),
\end{align*}
where $\Lambda=\{j\in \underline{p}\mid n_j=1\}$, $\Delta$ is a subset of $\Lambda$ of even cardinality and
\[P_{\Delta}(\alpha)=\prod_{j\in\Delta}\langle h_{i_j},\alpha\rangle.\]
\end{Theorem}

\begin{proof}
 Using Lemma \ref{vacuum traces} and \cite[Theorem 22]{Mason20231-PointAlgebras}, we have
 \begin{align*}
 Z_{V_L^+}(u,\tau)={}&\frac{1}{2}\sum_{\Delta\subseteq \Lambda}\frac{\theta_L(\tau,P_{\Delta})}{\eta(\tau)^k}\bigg(\sum_{\sigma\in \operatorname{Inv}_0(\underline{p}\backslash \Delta)}\prod_{(rs)}\delta_{i_r,i_s}\hat{E}_{n_r+n_s}(\tau)\bigg)\\
 &+\frac{1}{2}\eta(\tau)^{k/2}\left(\frac{\Theta_1(\tau)}{2}\right)^{-k/2}\bigg(\sum_{\sigma\in \operatorname{Inv}_0(\underline{p})} \prod_{(rs)}\delta_{i_r,i_s}\hat{F}_{n_r+n_s}(\tau) \bigg).
 \end{align*}
 Further, using Theorem \ref{At Heisenberg states}, Lemma \ref{vacuum traces}, Corollary~\ref{corollary1} and the above equation, we have the required result.
\end{proof}

\begin{Theorem}
Suppose $V$ is the $\mathbb{Z}_2$-orbifold of the lattice VOA $V_L$ corresponding to the $(-1)$-involution as defined in Section {\rm\ref{sec2.3}}, then for positive integers $n_i\geq 1$, $\alpha\in 2L$, corresponding to the state
$u\coloneqq h_{i_1}[-n_1]\cdots h_{i_p}[-n_p]e_{\alpha}\in V$,
where $e_{\alpha}$ is $f_{\alpha}$ or $g_{\alpha}$ accordingly as $p$ is even or odd, respectively,
\begin{align*}
Z_V(u,\tau)={}&\eta(\tau)^{k/2}\left(\frac{\Theta_1(\tau)}{2}\right)^{\langle\alpha,\alpha\rangle-k/2}\biggl(\sum_{\sigma} \prod_{(rs)(t)}\delta_{i_r,i_s}\Tilde{F}_{n_t}(\tau)\hat{F}_{n_r+n_s}(\tau)\biggr)\\
&+\eta(\tau)^{k/2}\left(\frac{\Theta_2(\tau)}{2}\right)^{\langle\alpha,\alpha\rangle-k/2}
\biggl(\sum_{\sigma}\prod_{(rs)(t)}\delta_{i_r,i_s}\Tilde{\overline{E}}_{n_t}\overline{E}_{n_r+n_s}(\tau)\biggr)\\
&+(-1)^l\eta(\tau)^{k/2}\left(\frac{\Theta_3(\tau)}{2}\right)^{\langle\alpha,\alpha\rangle-k/2}\biggl(\sum_{\sigma}\prod_{(rs)(t)}\delta_{i_r,i_s}
\Tilde{\overline{F}}_{n_t}\overline{F}_{n_r+n_s}(\tau)\biggr),
\end{align*}
where $\sigma$ ranges over all involutions of the index set $S_p=\{1,2,\dots, p\}$ and $(rs)$, $(t)$ range over the $2$-cycles and $1$-cycles respectively in the decomposition of $\sigma$ in $S_p$ as a product of disjoint $2$-cycles and $1$-cycles. For $\alpha\in L\backslash 2L$,
$u\coloneqq h_{i_1}[-n_1]\cdots h_{i_p}[-n_p]e_{\alpha}$,
we have
$
Z_V(u,\tau)=0$.
\end{Theorem}

\begin{proof}
 Using Lemma \ref{lattice element traces} and \cite[Theorem 19]{Mason20231-PointAlgebras}, we have
 \[Z_{V_L^+}(u,\tau)=\eta(\tau)^{k/2}\left(\frac{\Theta_1(\tau)}{2}\right)^{\langle\alpha,\alpha\rangle-k/2}\biggl(\sum_{\sigma} \prod_{(rs)(t)}\delta_{i_r,i_s}\Tilde{F}_{n_t}(\tau)\hat{F}_{n_r+n_s}(\tau)\biggr).\]
 Further, using Theorems \ref{Outside 2L} and \ref{Inside 2L}, Lemma \ref{lattice element traces}, we can compute $Z_{(V_L^T)^+}(u,\tau)$. Combining the two traces, we have the required result.
\end{proof}

\section[Modular invariance in the Z\_2-orbifold]{Modular Invariance in the $\boldsymbol{\mathbb{Z}_2}$-orbifold}\label{sec6}

Before we prove modular invariance in the orbifold, we shall state a lemma which is an immediate consequence of \cite[Proposition 23]{Mason20231-PointAlgebras}.

\begin{Lemma}\label{lemma G}
 For a state
$u=h_{i_1}[-n_1]\cdots h_{i_p}[-n_p]\mathbf{1}$,
 let \[G(u,\tau)=\sum_{\Delta\subseteq \Lambda}\frac{\theta_L(\tau,P_{\Delta})}{\eta(\tau)^k}\biggl(\sum_{\sigma\in \operatorname{Inv}_0(\underline{p}\backslash \Delta)}\prod_{(rs)}\delta_{i_r,i_s}\hat{E}_{n_r+n_s}(\tau)\biggr),\]
 where $\Lambda=\{ j\in \underline{p}\mid n_j=1\}$ and $\Delta$ is a subset of $\Lambda$ of even cardinality.
 When $L$ is unimodular of rank $k$, $G(u,\tau)$ is modular of weight equal to $L[0]$-weight of $u$.
\end{Lemma}

\begin{Proposition}\label{proposition1}
 Suppose $V$ is the $\mathbb{Z}_2$-orbifold of the lattice VOA $V_L$ corresponding to the $(-1)$-involution as defined in Section~{\rm\ref{sec2.3}}, then for every homogeneous state $u$ {\rm(}with respect to~$L[0]$-grading$)$ in $V$, we have
 \[Z_V(u,S\tau)=\tau^{wt(u)}Z_V(u,\tau),\qquad
 Z_V(u,T\tau)={\rm e}^{-2k\pi {\rm i}/24}Z_V(u,\tau),\]
 where $S,T\in {\rm SL}_2(\mathbb{Z})$ represent as usual the elements \smash{$\bigl(\begin{smallmatrix}
 0 & -1\\
 1 & 0
 \end{smallmatrix}\bigr)$} and
 \smash{$\bigl(\begin{smallmatrix}
 1 & 1\\
 0 & 1
 \end{smallmatrix}\bigr)$} respectively and hence $Z_V(u,\tau)$ is modular of weight $wt(u)$ {\rm(}with respect to $L[0])$ over the full modular group ${\rm SL}_2(\mathbb{Z})$ up to a~character.
\end{Proposition}

\begin{proof}
 Observe that the proposition is true for any homogeneous state $u$ in $\bigl(V_L^T\bigr)^+$ since for such states $Z_V(u,\tau)=0$. Thus, we prove the proposition for states in $V_L^+$ below.
 The $S$-invariance follows from the following:
 \begin{gather*}
 \hat{F}_k(S\tau)=\tau^k\overline{E}_k(\tau), \qquad \overline{E}_k(S\tau)=\tau^k\hat{F}_k(\tau), \qquad \overline{F}_k(S\tau)=\tau^k\overline{F}_k(\tau),\\
\Theta_1(S\tau)=(-{\rm i}\tau)^{1/2}\Theta_2(\tau), \qquad \Theta_2(S\tau)=(-{\rm i}\tau)^{1/2}\Theta_1(\tau), \qquad \Theta_3(S\tau)=(-{\rm i}\tau)^{1/2}\Theta_3(\tau),\\
\eta(S\tau)=(-{\rm i}\tau)^{1/2}\eta(\tau).\end{gather*}
Using Lemma \ref{lemma G}, we also have
$G(u,S\tau)=\tau^{wt(u)}G(u,\tau)$.
The $T$-invariance follows from the following:
 \begin{gather*}
 \hat{F}_k(T\tau)=\hat{F}_k(\tau), \qquad \overline{F}_k(T\tau)=\overline{E}_k(\tau),\qquad
 \hat{E}_k(T\tau)=\hat{E}_k(\tau), \qquad \overline{E}_k(T\tau)=\overline{F}_k(\tau),\\
 \Theta_1(T\tau)={\rm e}^{\frac{3\pi {\rm i}}{12}}\Theta_1(\tau), \qquad \Theta_2(T\tau)=\Theta_3(\tau), \qquad \Theta_3(T\tau)=\Theta_2(\tau),\\
 \eta(T\tau)={\rm e}^{\frac{\pi {\rm i}}{12}}\eta(\tau).\end{gather*}
 Using Lemma \ref{lemma G}, we also have
 \smash{$G(u,T\tau)={\rm e}^{-2k\pi {\rm i}/24}G(u,\tau)$}.
 Thus, for a Heisenberg state $u=h_{i_1}[-n_1]\cdots h_{i_p}[-n_p]\mathbf{1}$,
 \begin{align*}
 Z_V(u,T\tau)={}&\frac{1}{2}G(u,T\tau)+\frac{1}{2}\eta(T\tau)^{k/2}\left(\frac{\Theta_1(T\tau)}{2}\right)^{-k/2}\biggl(\sum_{\sigma\in \operatorname{Inv}_0(\underline{p})} \prod_{(rs)}\delta_{i_r,i_s}\hat{F}_{n_r+n_s}(T\tau) \biggr)\\
 &+\frac{1}{2}\eta(T\tau)^{k/2}\left(\frac{\Theta_2(T\tau)}{2}\right)^{-k/2}\biggl(\sum_{\sigma\in \operatorname{Inv}_0(\underline{p})}\prod_{(rs)}\delta_{i_r,i_s}\overline{E}_{n_r+n_s}(T\tau)\biggr)\\
 &+\frac{(-1)^l}{2}\eta(T\tau)^{k/2}\left(\frac{\Theta_3(T\tau)}{2}\right)^{-k/2}\biggl(\sum_{\sigma\in \operatorname{Inv}_0(\underline{p})}\prod_{(rs)}\delta_{i_r,i_s}\overline{F}_{n_r+n_s}(T\tau)\biggr)\\
 ={}&\frac{{\rm e}^{-2k\pi {\rm i}/24}}{2}G(u,\tau)+\frac{{\rm e}^{-2k\pi {\rm i}/24}}{2}\eta(\tau)^{k/2}\left(\frac{\Theta_1(\tau)}{2}\right)^{-k/2}\\
 &\phantom{+}{}\times \biggl(\sum_{\sigma\in \operatorname{Inv}_0(\underline{p})} \prod_{(rs)}\delta_{i_r,i_s}\hat{F}_{n_r+n_s}(\tau) \biggr)\\
 &+\frac{(-1)^l.(-1)^l.{\rm e}^{k\pi {\rm i}/24}}{2}\eta(\tau)^{k/2}\left(\frac{\Theta_3(\tau)}{2}\right)^{-k/2}\biggl(\sum_{\sigma\in \operatorname{Inv}_0(\underline{p})}\prod_{(rs)}\delta_{i_r,i_s}\overline{F}_{n_r+n_s}(\tau)\biggr)\\
 &+\frac{(-1)^l{\rm e}^{k\pi {\rm i}/24}}{2}\eta(\tau)^{k/2}\left(\frac{\Theta_2(\tau)}{2}\right)^{-k/2}\biggl(\sum_{\sigma\in \operatorname{Inv}_0(\underline{p})}\prod_{(rs)}\delta_{i_r,i_s}\overline{E}_{n_r+n_s}(\tau)\biggr).
 \end{align*}
 Rewriting $(-1)^l$ as ${\rm e}^{k\pi {\rm i}/8}$, we have
\smash{$Z_V(u,T\tau)={\rm e}^{-2k\pi {\rm i}/24}Z_V(u,\tau)$}.
Similarly, for a lattice state $u=h_{i_1}[-n_1]\cdots h_{i_p}[-n_p]e_{\alpha}$, we have the $T$-invariance given by
\[
Z_V(u,T\tau)={\rm e}^{-2k\pi {\rm i}/24}Z_V(u,\tau).
\tag*{\qed}
\]\renewcommand{\qed}{}
\end{proof}

\subsection*{Acknowledgments}

I would like to express my sincere gratitude to my advisor, Geoffrey Mason, for his invaluable guidance and support throughout this research. I also extend my heartfelt congratulations to Professor James Lepowsky on the occasion of his 80th birthday, and acknowledge his seminal contributions to the theory of vertex operator algebras. I would also like to thank the anonymous referees for their careful reading of the manuscript and for their helpful comments and suggestions, which significantly improved the presentation of this paper.

\pdfbookmark[1]{References}{ref}
\LastPageEnding

\end{document}